\documentclass[english,12pt]{article}
\usepackage{color}
\usepackage{verbatim}
\usepackage{float}
\usepackage{mathtools}
\usepackage{todonotes} 
\usepackage{amsmath,amssymb}
\usepackage{amsthm,enumerate}
\usepackage{hyperref}
\usepackage[capitalise, nameinlink]{cleveref}
\usepackage{amssymb}
\usepackage{fullpage}
\usepackage{xargs}[2008/03/08]

\colorlet{darkred}{red!80!black}
\colorlet{darkblue}{blue!80!black}
\colorlet{darkgreen}{green!60!black}

\usepackage[noend]{algpseudocode}
\floatstyle{ruled}
\newfloat{algorithm}{tbp}{loa}
\floatname{algorithm}{Procedure}

\usepackage{thmtools}

\theoremstyle{plain}
\newtheorem{thm}{Theorem}[section]
\newtheorem{cor}[thm]{Corollary}
\newtheorem{observation}[thm]{Observation}
\newtheorem{lem}[thm]{Lemma}
\newtheorem{conjecture}[thm]{Conjecture}
\newtheorem{prop}[thm]{Proposition}

\theoremstyle{remark}

\theoremstyle{definition}
\newtheorem{defn}[thm]{Definition}

\newtheorem{rem}[thm]{Remark}
\newtheorem{question}[thm]{Question}

\crefname{thm}{Theorem}{Theorems} 
\crefname{prop}{Proposition}{Propositions} 
\crefname{conjecture}{Conjecture}{Conjectures} 
\crefname{observation}{Observation}{Observations} 
\crefname{cor}{Corollary}{Corollaries} 
\crefname{claim}{Claim}{Claims} 
\crefname{lem}{Lemma}{Lemmas} 
\crefname{rem}{Remark}{Remarks}

\newcommand{\ourtitle}{Non-induced inherent graphs}

\colorlet{myGreen}{green!70!black}
\definecolor{myBlue}{rgb}{0.25, 0.0, 1.0}
\definecolor{lgray}{rgb}{0.75, 0.75, 0.75}
\hypersetup{
  colorlinks=true,
  linkcolor=myBlue,
  citecolor=myGreen,
  urlcolor=myGreen,
  bookmarksopen=true,
  bookmarksnumbered,
  bookmarksopenlevel=2,
  bookmarksdepth=3
  pdftitle = {\ourtitle},
  pdfauthor= {},
}

\newcommand{\cpi}[1]{\langle #1 \rangle}
\newcommand{\cp}[3]{(\cpi {#1},\cpi {#2},\cpi {#3})}

\newcommandx\s[1][usedefault, addprefix=\global, 1=]{\underset{#1}{\rightarrow}}%

\newcommand{\mpfile}[2]{\includegraphics{#1-#2-mps.pdf}}
\let\es\varnothing

\def\cH{\mathcal H}
\def\cF{\mathcal F}
\let\ph\varphi
\let\sm\setminus
\usepackage{babel}
\title{\ourtitle\thanks{This work is supported in part by the Slovenian Research and Innovation Agency (I0-0035, research program P1-0285 and research projects J1-60012, J1-70035, J1-70046, and N1-0370). The first and the fourth authors were supported by a research project HSE-BR-2025-024 implemented as part of the Basic Research Program at HSE University.}}
\author{
Vladimir Gurvich\thanks{National Research University Higher School of Economics and RUTCOR, Rutgers University, New Jersey, USA, \texttt{vladimir.gurvich@gmail.com}.}\and
Matjaž Krnc\thanks{FAMNIT and IAM, University of Primorska, \texttt{matjaz.krnc@upr.si}.}\and
Martin Milanič\thanks{FAMNIT and IAM, University of Primorska, \texttt{martin.milanic@upr.si}.}\and
Mikhail Vyalyi\thanks{Federal Research Center “Computer Science and Control” of the Russian Academy of Science; National Research University Higher School of Economics; Moscow Institute of Physics and Technology;  \texttt{vyalyi@gmail.com}.}}
\begin{document}
\maketitle

\begin{abstract}
Concepts  of  avoidable paths and inherent two-rooted graphs were recently introduced for the case of induced subgraphs.
In this paper, we translate this framework to the setting of not necessarily induced (NNI) subgraphs, where the structural constraints, as well as conclusions, are significantly relaxed. 
We develop a new sufficient condition for inherence, and use it to provide complete characterizations of inherent two-rooted paths and forks. 
We show that for these families, the obvious necessary condition for inherence from the induced case is sufficient in the NNI setting. 
To complement these results, we introduce two general methods for proving non-inherence and apply them to demonstrate the non-inherence of several classes of two-rooted trees.

\bigskip


\bigskip

\end{abstract}

\section{Introduction}
In graph theory, graph properties are often characterized by forbidden or mandatory substructures. These substructures are analyzed in two main ways: as induced subgraphs, which preserve all adjacencies from the parent graph, or as the more flexible \emph{not necessarily induced} (NNI) subgraphs, which require only a subset of vertices and edges.

The concept of ``avoidability," central to our work, has its roots in the study of chordal graphs. 
The foundational result is Dirac's theorem \cite{MR130190}, which states that every chordal graph contains a simplicial vertex—a vertex whose neighborhood forms a clique. 
Recognizing that not all graphs have simplicial vertices, Ohtsuki, Cheung, and Fujisawa introduced in 1976 the concept of an \textsl{avoidable vertex}: one for which every induced path of length two centered at this vertex is part of an induced cycle~\cite{MR485552}.
They proved that such a vertex exists in every graph.
The term `avoidable vertex' was introduced by Beisegel et al.~\cite{BCGMS19}; an alternative term used in the literature is `OCF-vertex' (see, e.g.,~\cite{MR2057267,MR2191642,MR2204109}).

Nearly half a century later, this notion was generalized from vertices to paths. 
Beisegel et al.~\cite{BCGMS19,BCGMS19DAM} defined an avoidable path as an induced path whose every ``extension'' -- adding a new vertex to each endpoint so that the obtained longer path is still induced -- is contained in an induced cycle.  
They conjectured that every graph containing an induced $k$-vertex path must also contain an avoidable $k$-vertex path and proved the claim for $k = 2$.
The conjecture of Beisegel et al.~generalizes not only the aforementioned result of Ohtsuki, Cheung, and Fujisawa~\cite{MR485552} but also a result by Chv\'atal, Rusu, and Sritharan~\cite{MR1927566} who proved the conjecture for graphs excluding chordless cycle of length at least~$k+3$.
In full generality, the conjecture was proved by Bonamy et al.~\cite{MR4245221} and later strengthened by Gurvich et al.~\cite{GKMV20+}, who showed that every induced path can be transformed into an avoidable one through a series of local transformations called shifts. 

A further generalization was provided in \cite{gurvich2023avoidability}, which established the framework of inherent two-rooted graphs for the induced case. 
A \emph{two-rooted graph} $(H,s,t)$ consists of a graph $H$ with two specified vertices $s$ and $t$, called \emph{roots}.  

This framework formalizes the notion of an extension for arbitrary structures: an \emph{extension} of an induced copy of $(H,s,t)$ in a larger graph $G$ is an induced subgraph obtained by adding to the copy of $H$ a pair of pendant edges incident with the copies of $s$ and $t$, respectively.
Such an extension is \emph{closable} if a path exists in $G$ between these new pendant vertices, using only vertices outside the closed neighborhood of the original copy. 
A copy of $(H,s,t)$ is then defined as \emph{avoidable} if all of its possible extensions are closable.
In particular, a copy is avoidable if it admits no extensions, in which case it is called \emph{simplicial}.
Finally, a two-rooted graph is \emph{inherent} if any graph containing an induced copy of $H$ must also contain an avoidable induced copy of $(H,s,t)$.\footnote{For more precise definitions see~\Cref{subsec:definitions}.}
Gurvich et al.~\cite{gurvich2023avoidability} identified an important necessary condition for a connected two-rooted graph to be inherent (\emph{NCI} for short).

\begin{defn}\label{def:NCI}
A two-rooted graph $( H, s, t)$ is said to have \emph{NCI property} if:
\begin{itemize}
    \item either $H$ is a path with an endpoint $ s$, 
    \item or $H$ is a tree with maximum degree~$3$ such that $1 = d_{ H}( s) \le d_{ H}( t)\le 2$, with $d_{ H}( t) = 1$ if and only if $ s =  t$.
\end{itemize}
\end{defn}

In particular, a two-rooted graph $(H, s, t)$ can only have the NCI property if it is a \emph{two-rooted tree}, that is, $H$ is a tree. 
In \cite[Theorem 1.8]{gurvich2023avoidability} it was shown that every inherent connected graph has the NCI property. 
Yet, determining inherence for a specific graph seems to be difficult. 
So far, the only known connected (induced) inherent graphs are paths rooted in their endpoints (see \cite{MR4245221}), and even some small examples, such as the two-edge path with adjacent roots, are not yet settled (see~\cite{gurvich2023avoidability}). 

Furthermore, it is an open question whether the problem of  recognizing if a given two-rooted graph is inherent
is decidable.

\subsection{Our results}

Recall that NNI stands for ``not necessarily induced''. 
The present work pivots from the induced to the NNI setting. 
This transition is a natural and crucial step, as the structural constraints are significantly relaxed, in the sense that a non-induced copy of a graph $H$ can exist in many more contexts than an induced one; 
on the other hand, in general, there are more extensions in the NNI setting than in the induced one.
We adopt the two-rooted graph framework but redefine the core concepts for NNI subgraphs. 
Our goal is to determine how this added flexibility impacts the property of inherence and to identify classes of two-rooted graphs that are robust enough to be inherent in this more general setting.

This paper provides the first systematic study of inherence in the NNI setting. 
Our starting point is the observation that adapting the definitions and the proofs of \hbox{\cite[Section 3]{gurvich2023avoidability}} \textit{(mutatis mutandis)} to the NNI setting shows similar results also in the NNI setting.
In particular, the following holds. 

\begin{observation}\label{obs:NCI-NNI}
    All connected inherent two-rooted graphs have NCI property.
\end{observation}

We develop new proof techniques tailored to the NNI setting and establish several foundational results. 
Our main contributions are as follows.

A powerful sufficient condition for a copy of $(H,s,t)$ to be avoidable in a graph $G$ is that the subgraph induced by the vertices not in the copy is connected. 
Any extension of such a copy would have its new vertices within this subgraph and would therefore be closable; in particular, the copy is avoidable. 
It turns out that for certain two-rooted graphs, the mere absence of simplicial copies guarantees the existence of such a copy of $H$. 
Such two-rooted graphs are called \emph{complementary connected or simplicial} (CCS) graphs. 
By definition, CCS graphs are inherent. 
The converse implication is open.

Using the CCS framework, we provide complete characterizations of inherence for all two-rooted paths\footnote{A \emph{two-rooted path} is a two-rooted tree $(H,s,t)$ such that $H$ is a path.} and for all two-rooted trees with at most six vertices.
\begin{description}
    \item[Paths:] We call a path with roots in both endpoints an \emph{endpoint-rooted path}. 
    Gurvich et al.\ showed the following result in the NNI setting.
\begin{thm}[\cite{GKMV20+}]\label{thm:endpoints}
  Any endpoint-rooted path is inherent.
\end{thm}
We show that for two-rooted paths, the necessary condition for inherence (NCI) is also sufficient in the NNI setting.
\begin{thm}\label{thm:paths_scc}
Every two-rooted path satisfying the NCI property is a CCS graph, and therefore inherent.
\end{thm}
\item[\boldmath Small trees:]
We completely characterize inherent two-rooted trees with up to $6$ vertices.
It turns out that, when the underlying tree is a \emph{fork}, that is, the unique $5$-vertex tree with maximum degree $3$, inherence is equivalent to the CCS property.
Furthermore, apart from paths and forks,\footnote{A \emph{two-rooted fork} is a two-rooted tree $(H,s,t)$ such that $H$ is a fork.} the only inherent two-rooted tree on at most $6$ vertices is $T_1(2,1,2)$, depicted in \cref{pic:T1212-names}(a) and formally defined in \cref{subsec:definitions}. 
\begin{restatable}{thm}{thmsixvertex}\label{thm:T1212}
    $T_1(2,1,2)$ is inherent.
\end{restatable}
\end{description}

As mentioned above, apart from paths and forks, among connected graphs on at most~$6$ vertices there exists a single inherent one. 
Indeed, it seems that inherent graphs are very rare. 
In this regard, we provide additional necessary conditions for inherence. 
We show that a two-rooted tree cannot be inherent if its bipartition is too unbalanced. 
Furthermore, we present a \textsl{path augmentation} method to construct confining graphs for several families of two-rooted trees, including combs of type II, rakes, and certain two-rooted leaf-extended trees (see \cref{subsec:definitions,sec:direct-confinement} for definitions), thereby proving their non-inherence in the NNI setting.

\subsection{Preliminaries}
\label{subsec:definitions}

While we mostly work with undirected graphs, we also use \emph{arcs} for a particular orientation of a given edge.
A \emph{copy} of a graph $H$ in a graph $G$ is a subgraph of $G$ isomorphic to $H$.
A \emph{subgraph embedding} of a graph $H$ in a graph $G$ is an injection $\ph\colon V(H)\to V(G)$ that maps adjacent vertices to adjacent vertices. 
The path graph on $n$ vertices is denoted by $P_n$.
For other standard graph theory notation we refer the reader to  Diestel~\cite{diestel-graph-theory}.

A \emph{two-rooted graph}  is a triple $(H, s, t)$ such that $H$ is a graph and $ s$ and $t$ are two (not necessarily distinct) vertices of $H$ which we call \emph{roots}. We will denote two-rooted graphs in calligraphic font, such as $\mathcal{G,H}$. 
For convenience and without loss of generality, we always assume for a two-rooted graph $(H, s, t)$ that $d_{H}(s)\leq d_{H}(t)$.

In what follows we adapt the definitions of terms from~\cite{gurvich2023avoidability} such as \emph{copy, extension, avoidability, inherence}, and \emph{confining graph} to the NNI setting. 
We point out that in the induced setting these notions differ from ours, see \cref{rem:terminology}.
Given a graph $G$ and a two-rooted graph $( H,  s,  t)$, a \emph{copy of $( H,  s,  t)$ in $G$} is any two-rooted graph $\cp Hst$ such that $\cpi H$ is a subgraph of~$G$ for which there exists an isomorphism of $H$ to $\cpi H$ mapping $s$ to $\cpi s$ and $t$ to $\cpi t$.

Given a two-rooted graph $\cH =(H,s,t)$, the \emph{extension graph} of $\cH$ is the graph $H'$  obtained from $H$ by adding to it two pendant edges $ss'$ and $tt'$; the corresponding two-rooted graph $\cH' = (H',s',t')$ is said to be an \emph{extension of $\cH$}.
Given a graph $G$, a two-rooted graph $( H,  s,  t)$, and a copy $\cp Hst$ of it in $G$, an \emph{extension of $\cp Hst$ in $G$} is any extension $(H',s',t')$ of $\cp Hst$ such that $H'$ is a subgraph of $G$.
In particular, $V(H') = V(\cpi H)\cup\{s',t'\}$, vertices $s'$ and $t'$ are distinct, the graph obtained from $H'$ by deleting $s'$ and $t'$ is $\cpi H$, and $\cpi s$ and $\cpi t$ are unique neighbors of $s'$ and $t'$ in $H'$, respectively. 
Furthermore, we say that an extension $(H',s',t')$ of $\cp Hst$ in $G$ is \emph{closable} if there exists an  $s',t'$-path in $G$ having no vertex in common with $V(\cpi H)$. 

\begin{defn}\label{def:avoidable}
Let $\cp Hst$ be a copy of a two-rooted graph $( H,  s,  t)$ in a graph $G$.
We say that $\cp Hst$ is \emph{avoidable} (in $G$) if all extensions of $\cp Hst$ in $G$ are closable.
\end{defn}

In particular, a copy $\cp Hst$ of $( H,  s,  t)$ in $G$ that has no extensions is trivially avoidable.
Such a copy $\cp Hst$ will be called {\em simplicial}.

\begin{rem}
If $H$ is a path and $ s$ and $ t$ are its endpoints, the above definitions of an extension of a copy  of $( H, s,  t)$ in a graph $G$, a closable extension, a simplicial copy, and an avoidable copy coincide with the corresponding definitions from Gurvich et al.~\cite[Section 1.1]{GKMV20+} for the case $t = \mathtt{pth}$.
\end{rem}

Finally, we introduce the core definition of this paper. In fact, this is the definition from~\cite{gurvich2023avoidability} adjusted to the NNI setting.

\begin{defn}\label{inherent}
A two-rooted graph $( H,  s,  t)$ is \emph{inherent} if every graph $G$ that contains a copy of $H$ also contains an avoidable copy of $(H,s,t)$.
\end{defn}

The following concept provides a natural certificate of non-inherence.

\begin{defn}\label{def-confining-graph}
Given a two-rooted graph $( H, s, t)$ and a graph $G$, we say that $G$ \emph{confines} $( H, s, t)$ (or: is a \emph{confining graph} for $( H, s, t)$) if $G$ contains a copy of $ H$ but no avoidable copy of $( H, s, t)$.
\end{defn}

In~\cite{gurvich2023avoidability} a necessary condition for inherence of two-rooted graphs was introduced, called \emph{PE-inherence}. 
The condition remains necessary in the NNI setting.
It follows from~\cite[Observation 3.24 and Proposition 3.26]{gurvich2023avoidability} that a two-rooted tree is PE-inherent if and only if it cannot be confined by a tree.\footnote{Necessity of the condition follows from~\cite[Observation 3.24]{gurvich2023avoidability}, while sufficiency follows from a generalization of~\cite[Proposition 3.26]{gurvich2023avoidability} stating that for every two-rooted tree $(H,s,t)$, every confining tree of $(H,s,t)$ contains a subgraph isomorphic to the limit graph of $(H,s,t)$.
This can be proved with the obvious adaptation of the proof of~\cite[Proposition 3.26]{gurvich2023avoidability} and implies that no PE-inherent two-rooted tree can be confined by a tree. See \cite{gurvich2023avoidability} for definitions.}

Let us recall the definitions of several infinite families of PE-inherent two-rooted trees introduced in~\cite{gurvich2023avoidability}. 
A graph $ F(l,k,l')$ consists of a path 
$P_{l+k+l'}=(a_1,\dots,a_{l+k+l'})$ and $k$ pending edges added to vertices ${a_{l+1},\dots,a_{l+k}}$.

\begin{itemize}
\item \emph{Two-rooted combs of type I} are the two-rooted graphs $T_1(l,k,l') = (F(l,k,l'), a_1, a_l)$, where $l\ge 1$ and $k,l'\ge 0$ (see \cref{pic:T1-211,pic:T1-112} for two examples, namely $T_1(2,1,1)$ and $T_1(1,1,2)$, respectively).
Note that $T_1(l,0,l')$ is a two-rooted path and $T_1(l,0,0)$ is an endpoint-rooted path.
\item 
\emph{Two-rooted combs of type II} are the two-rooted graphs $T_2(p,q)=(F(p,q,p),a_1,a_{p+q+1})$, $p,q\ge 1$ (see \cref{fig:combII-confined} for an example, namely $T_2(3,3)$).
\item 
\emph{Two-rooted combs of type  III} are the two-rooted graphs  
$$T_3(p,q) = (F(p,q,p+1),a_1,a_{p+q+1}),$$ $p\ge 1$ and $q\ge 0$ (see
\cref{pic:T3-11} for an example, namely $T_3(1,1)$).
\item For an integer $q\ge 2$ we denote by $T_4(q)$ the two-rooted graph obtained from the two-rooted comb $T_2(2,q)$ of type II by subdividing each edge of the form $a_ib_i$, $i = 3,\ldots, q+2$. Any two-rooted graph equivalent to $T_4(q)$, for some $q$, will be referred to as a \emph{two-rooted rake} (see~\cref{fig:rakes-confined}(a) for an example, namely $T_4(3)$).
\item For an integer $q\ge 3$ we denote by $T_5(q)$ the two-rooted graph obtained from the two-rooted comb $T_2(2,q)$ of type II by subdividing each edge of the form $a_ib_i$, $i= 3,\ldots, q+2$, and adding a pendant edge to each vertex of degree two joining $a_i$ and $b_i$ for all $i \in \{4,\ldots, q+1\}$.
Any two-rooted graph equivalent to $T_5(q)$, for some $q$, will be referred to as a \emph{two-rooted split rake} (see~\cref{fig:rakes-confined}(b) for an example, namely $T_5(4)$).
\end{itemize}

\subsection{Structure of the paper}

The paper is structured as follows.
In \cref{sec:inherent} we develop a framework for proving inherence (see \Cref{fig:implications} for a summary).
We present proofs of \cref{thm:paths_scc}, the results on two-rooted forks, and \cref{thm:T1212} in \cref{sec:paths,sec:fork,sec:T1212}, respectively, while in \cref{sec:disconn} we discuss disconnected inherent graphs.
\cref{sec:confining} presents two main techniques for constructing confining graphs, which are used to prove that a two-rooted graph is not inherent. 
Namely, \cref{subsec:conf-bipartite} demonstrates how to use complete bipartite graphs to confine trees with imbalanced color classes, while 
\cref{sec:direct-confinement} introduces the so-called \textsl{path augmentation} method. 
This technique is then used to construct confining graphs for several families, including combs of type II, rakes, and certain two-rooted leaf-extended trees.
In \cref{sec:conclusion} we outline possible future research directions and pose some open questions.
In Appendix~\ref{sec:dictionary}, we collect the abbreviations of various properties of two-rooted graphs defined in \cref{sec:inherent}, while in Appendix~\ref{sec:small-graphs}, we provide the analysis of all non-inherent two-rooted graphs of order at most~$6$.

\bigskip
\noindent{\bf AI Usage Statement:} Artificial Intelligence was not used in any way for this research or writing. 

\section{A framework for proving inherence}\label{sec:inherent}

A~simplicial copy of a two-rooted graph is avoidable by definition. 
So, in the pursuit of inherent graphs we restrict our attention to  graphs that contain copies of a two-rooted graph but all copies are non-simplicial. 
Also, we are mainly interested in connected two-rooted graphs. 
Thus we proceed by defining several concepts that capture these properties, and describe relations between them.

\begin{defn}\label{def:H-smooth}
  Let $\cH =(H,s,t)$ be a two-rooted graph. 
  A graph $G$ is \emph{$\cH$-smooth} if $G$ is connected, $G$ contains a copy of $H$, and there are no simplicial copies of  $\cH$ in $G$.
 \end{defn}

Due to \cref{obs:NCI-NNI} we are mainly concerned in $\cH$-smooth graphs for the case when $\cH$ has the NCI property.

In many cases $\cH$-smoothness implies that there are many copies of $\cH$. We introduce several concepts capturing this `richness' of copies of a two-rooted graph.
Recall that a \emph{subgraph embedding} of a graph $H$ in a graph $G$ is an injection $\ph\colon V(H)\to V(G)$ that maps adjacent vertices to adjacent vertices.

\begin{defn}\label{def:H-universality}
  A two-rooted graph $\cH =(H,s,t)$ is said to be
  \begin{itemize}
  \item \emph{VU} (Vertex Universal) if for any  $\cH$-smooth graph $G=(V,E)$, any $v\in V$, and any $x\in V(H)$, there exists a subgraph embedding  $\ph$ of $H$ in $G$ such that $\ph(x)=v$; 
  
  \item \emph{WVU} (Weakly Vertex Universal) if for any  $\cH$-smooth graph $G=(V,E)$, any $v\in V$ contained in a copy of $H$ in $G$, and any $x\in V(H)$, there exists a subgraph embedding  $\ph$ of $H$ in $G$ such that $\ph(x)=v$; 
  
\item \emph{AU} (Arc Universal) if for any  $\cH$-smooth graph $G=(V, E)$, any arc $(u,v)\in E$, and any arc $(x,y)\in E(H)$,  there exists a subgraph embedding $\ph$ of $H$ in $G$ such that $\ph(x)=u$ and $\ph(y) =v$; 

 \item \emph{WAU} (Weakly Arc Universal) if for any  $\cH$-smooth graph $G=(V, E)$, any arc $(u,v)\in E$ contained in a copy of $H$ in $G$, and any arc  $(x,y)\in E(H)$, there exists a subgraph embedding  $\ph$ of $H$ in $G$ such that  $\ph(x)=u$ and $\ph(y) =v$;  

\item \emph{AS} (Arc Separating) if for any  $\cH$-smooth graph $G=(V, E)$ and any arc $e=(x,y)\in E$ such that $e$ is contained in a copy of $H$ in $G$  there exists a copy of $H$ that contains $x$ and does not contain $y$.
  \end{itemize}
\end{defn}

For a VU graph $\cH$,
each vertex $v$ in an $\cH$-smooth graph $G$ is contained in a copy of $\cH$ such that $v$ corresponds to a vertex of the maximal degree in $H$. So, the following corollary is immediate. 

\begin{cor}\label{cor:deg3}
Let $\cH= (H,s,t)$ be VU and let $G$ be $\cH$-smooth. Then $\Delta(H)\le\delta(G)$.
\end{cor}

In general, the relations between pairs of concepts among (W)VU, (W)AU, and AS are open except for the trivial implications:  
VU implies WVU, AU implies WAU, AU implies VU, and WAU implies WVU.  

\begin{lem}\label{lm:H-cover}
For a two-rooted graph $\cH =(H,s,t)$ with $\delta(H) = 1$, WVU implies VU.
\end{lem}

\begin{proof}
  Let $V_0$ be the union of vertex sets of all copies of $H$ in an $\cH$-smooth graph $G=(V,E)$. 
  Then, $V_0\ne \es$.
  Let $xy\in E$ be an edge such that $x\in V_0$. 
  Since $\cH$ is WVU, there exists a subgraph embedding  $\ph\colon V(H)\to V(G)$ such that $\ph^{-1}(x)$ is  
  a neighbor of a pendant vertex $h$ in $H$. 
  If $y$ is contained in $\ph(V(H))$, then $y\in V_0$. 
  Otherwise, substituting $\ph(h)$ by $y$ we get another copy of $H$ in $G$ and this copy contains $y$. 
  Therefore $y\in V_0$. 
  So, $y\in V_0$ in any case. 
  Thus, $V_0=V$ since $G$ is connected. 
  This implies that $\cH$ is VU.
\end{proof}

In the inherence proofs below we need AS and we derive it from WAU. 
But we are not aware of a general implication `WAU implies AS', even for PE-inherent two-rooted trees.  
Nonetheless, the following partial result will be helpful. 
For the statement we need one more definition.

\begin{defn}\label{def:EM}
  A two-rooted graph $\cH=(H,s,t)$ is \emph{EM} (Extension Movable) if the extension graph of $\cH$ contains more than one copy of $H$. 
\end{defn}

It is clear from the definition that every PE-inherent two-rooted tree is EM (otherwise it would be confined by its extension graph).  

\begin{lem}\label{lm:EMandWAU->AS}
  If a connected two-rooted graph $\cH =(H,s,t)$ is EM and WAU, then $\cH$ is AS.
\end{lem}

\begin{proof}
  Let $G=(V,E)$ be a $\cH$-smooth graph and let $e=(x,y)$ be an arc of $G$ that is contained in a copy of $H$ in $G$. 
  Since $\cH$ is EM, in the extension graph $T$ of $\cH$ we can choose a copy $H'$ of $H$ different from $H$.
  Furthermore, we can choose an arc $(h_1,h_2)$ of $H$ such that $h_1\in V(H')$ and $h_2\notin V(H')$. 
  Such an arc does exist, since $T$ is connected, any edge in $E(H')\sm E(H)$ is pendant, thus, any path in $T$ from a vertex in $V(H')\setminus V(H)$ to a vertex in $V(H)\setminus V(H')$ has an arc that $(h_1,h_2)$ of $H$ such that $h_1\in V(H')$ and $h_2\notin V(H')$. 
  
  Since $\cH$ is WAU, there exists a copy $\cH'$ of $\cH$ in $G$ such that $h_1$ corresponds to $x$ and $h_2$ corresponds to $y$. 
  Since $G$ is $\cH$-smooth, $\cH'$ has an extension in $G$.
  Choose in the extension graph a copy $H''$ of $H$ such that $x\in V(H'')$ and $y\notin V(H'')$.
\end{proof}

In verifying WAU for a two-rooted graph $\cH =(H,s,t)$, it is helpful to use the following equivalence relation on the set of arcs of $H$. 

\begin{defn}\label{def:G-equiv}
For graphs $G$ and $H$, two arcs $e_1$ and $e_2$ of $H$ are \emph{$G$-equivalent} if there are two subgraph embeddings $\ph_1$, $\ph_2$ of $H$ in $G$ such that $\ph_1(e_1) = \ph_2(e_2)$. 
We denote this relation as $e_1\sim_G e_2$. 
\end{defn}

For any two graphs $G$ and $H$, the relation $\sim_G$ is an equivalence relation on the set of arcs of $H$. 
Indeed, it is reflexive and symmetric by definition. 
Furthermore, since a composition of isomorphisms is an isomorphism, this relation is transitive.\footnote{More specifically, if there exist subgraph embeddings $\ph_1$ and $\ph_2$ of $H$ in $G$ such that $\ph_1(e_1) = \ph_2(e_2)$, as well as subgraph embeddings $\psi_2$ and $\psi_3$ of $H$ in $G$ such that $\psi_2(e_2) = \psi_3(e_3)$, then $\psi_2\circ \ph^{-1}_2\circ \ph_1$ is a subgraph embedding of $H$ in $G$ such that $(\psi_2\circ \ph^{-1}_2\circ \ph_1) (e_1) = \psi_3(e_3)$.}

The following lemma follows immediately from the definitions. 
We state it explicitly for future reference. 

\begin{lem}\label{lm:equivalence}
  A two-rooted graph $\cH = (H,s,t)$ is WAU if for every $\cH$-smooth graph $G$, every two arcs of $H$ are $G$-equivalent to each other.
\end{lem}

Now we turn to an important property sufficient for the existence of an avoidable copy.
Let $\cp Hst$ be a copy of $(H,s,t)$ in a graph $G=(V,E)$ such that $G-V(\cpi H)$ is connected.
Then the copy is avoidable: every extension has endpoints in $V\setminus V(\cpi H)$ and they are connected by a path inside $V\setminus V(\cpi H)$.
Note that this property does not depend on roots.

\begin{defn}\label{def:HCC}
  Given two graphs $G=(V,E)$ and $H$, the graph $G$ has the $H$-\emph{Connected Complement} property) ($H$-\emph{CC property} for short) if there exists a copy $H'$ of $H$ in $G$ such that $G-V(H')$ is connected.
\end{defn}

\begin{defn}\label{def:CCS}
  A two-rooted graph  $\cH = (H,s,t)$ is a \emph{CCS graph} (\emph{Connected Complement or Simplicial}) if every $\cH$-smooth  graph has the $H$-CC property.
\end{defn}

It follows immediately from the definition that every CCS two-rooted graph is inherent.
More concepts and notations are needed in the proofs of the CCS property.

\begin{defn}
Let $H$ be a  graph.  
For a graph $G=(V,E)$, a pair $(U,R)$
of subsets of $V$ is  \emph{$H$-avoiding} if $\{U,R\}$ is a partition of $V$, the subgraph of $G$ induced by $U$ is connected, and the subgraph of $G$ induced by $R$ contains a copy of $H$.
An $H$-avoiding pair $(U,R)$  is called \emph{inclusion-maximal} if no $H$-avoiding pair $(U',R')$ satisfies $U\subset U'$.
\end{defn}

For graphs $G$ and $H$, an $H$-avoiding pair $(U,R)$ in $G$, and $x\in R$, we denote by $R_{x}$ the graph $G[R\setminus\{x\}]$.\footnote{For $S\subseteq V(G)$, we denote by $G[S]$ the subgraph of $G$ induced by $S$.}
Note that an inclusion-maximal $H$-avoiding pair in a graph $G $ satisfies the \emph{hitting condition}: if $x\in N_G(U)$ then $R_x$ contains no copies of $H$. 
In other words, each copy of $H$ in $G[R]$  contains $N_G(U)$. 
Otherwise, if $r\in  N_G(U)$ is not covered by the copy, then the set $U$ can be extended by a vertex $r$ preserving the properties of an $H$-avoiding pair. 
Later we will need the following property, which is a simple implication of the hitting condition.

\begin{prop}\label{URx-connected}
Let $\cH =(H,s,t)$ be a two-rooted tree, $(U,R)$ an inclusion-maximal $H$-avoiding pair in a $\cH$-smooth graph $G=(V,E)$, and $x\in N_G(U)$. 
Then $x$ has a neighbor in each component of $R_x$.
\end{prop}

\begin{proof}
Suppose for contradiction that ${N_G(x)\cap C = \es}$ for some component $C$ of $R_x$.\footnote{For simplicity, we do not distinguish between a connected component and its vertex set.} 
Since $G$ is connected, $C\cap N_G(U)\ne \es $.
Let $y\in C\cap N_G(U)$.
The hitting condition implies that there exists a copy $\cpi H$ of $H$ containing $x$ and $y$.
Let $P$ be a shortest path from $x$ to $y$ in $\cpi H$.
Then, all the vertices of $P$ except $x$ belong to $C$, implying that $N_G(x)\cap C \neq \es$, a~contradiction.
\end{proof}

\begin{defn}\label{def:HCS}
  A two-rooted graph $\cH =(H,s,t)$ is \emph{HCS} (H-Component is Small) if  for any $\cH$-smooth graph $G$ and any inclusion-maximal $H$-avoiding pair $(U,R)$ in $G$, the following implication holds: 
  if there exists a vertex $x\in N_G(U)$ such that $R_x$ is connected, then $|R| = |V(H)|$. 
\end{defn}

The HCS property implies the $H$-CC property of any $\cH$-smooth graph $G$, since by definition the graph $G[R]$ contains a copy of $H$.
The basic lemma for the proofs of inherence is as follows.

\begin{lem}\label{lm:AS-AU-HCS->CCS}
  If a two-rooted graph $\cH =(H,s,t)$ is AS, AU, and HCS, then it is CCS.
\end{lem}

As an immediate corollary we get that any AS, AU, and HCS graph is inherent.

\begin{proof}
  Take an $\cH$-smooth graph $G=(V,E)$ and an inclusion-maximal $H$-avoiding pair $(U,R)$ in $G$.
  
  We claim that for every vertex $x\in N_G(U)$, every connected component $C$ of $R_x$ satisfies $N_G(U)\cap C\ne \es$.
  We prove the claim by contradiction. 
  Suppose that $N_G(U)\cap C = \es$ for some  component $C$ of $R_x$, where $x\in N_G(U)$.   
  Then, $C$ is a component in $G-x$, since $N_G(U)\cap C= \es$. 
  It follows from \cref{URx-connected} that there exists an edge $xv\in E$ such that $v\in C$.  
  Since $\cH$ is AU, there exists a copy of $H$ in $G$ that contains the arc $(v,x)$.
  Hence, as $\cH$ is AS, there is a copy of $H$ that contains $v$ and does not contain $x$. 
  This copy should be contained in $C$ since $N_G(C)=\{x\}$.
  Thus $(U\cup \{x\}, R\setminus\{x\})$ is an $H$-avoiding pair, contradicting the maximality of $(U,R)$.

  For $x\in  N_G(U)$, let $R_{x,1}, R_{x,2}, \dots, R_{x,s}$ be the connected components of $R_x$ ordered non-increasingly by size, and let $r_{x,i}=|R_{x,i}|$ for all ${1\le i\le s}$.
Then $r_{x,1} \geq r_{x,2} \geq \ldots \geq r_{x,s}$.
Choose $x\in N_G(U)$ such that the vector $r_x= (r_{x,1}, r_{x,2}, \dots, x_{r,s})$ is maximal with respect to the lexicographic order. 
We claim that $R_x$ is connected. 
Since $\cH$ is HCS, this implies that $H$-CC property of $G$ does hold and, consequently, that $\cH$ is CCS.

Now, we prove the claim. Suppose for contradiction that $s\geq2$. 
By \cref{URx-connected}, $x$~has a neighbor $w$ in $R_{x,2}$.
Then, ${R_{x,2}\cap N_G(U) = \es}$ due to the maximality in the lexicographic order, since if $v\in R_{x,2}\cap N_G(U)$, then the vector $r_{v}$ is lexicographically larger than $r_x$ (the component $R_{x,1}$ is enlarged by at least  $x$ by \cref{URx-connected}). 
As we have seen before, this is impossible. 
Hence, we come to a contradiction.
\end{proof}

\begin{figure}
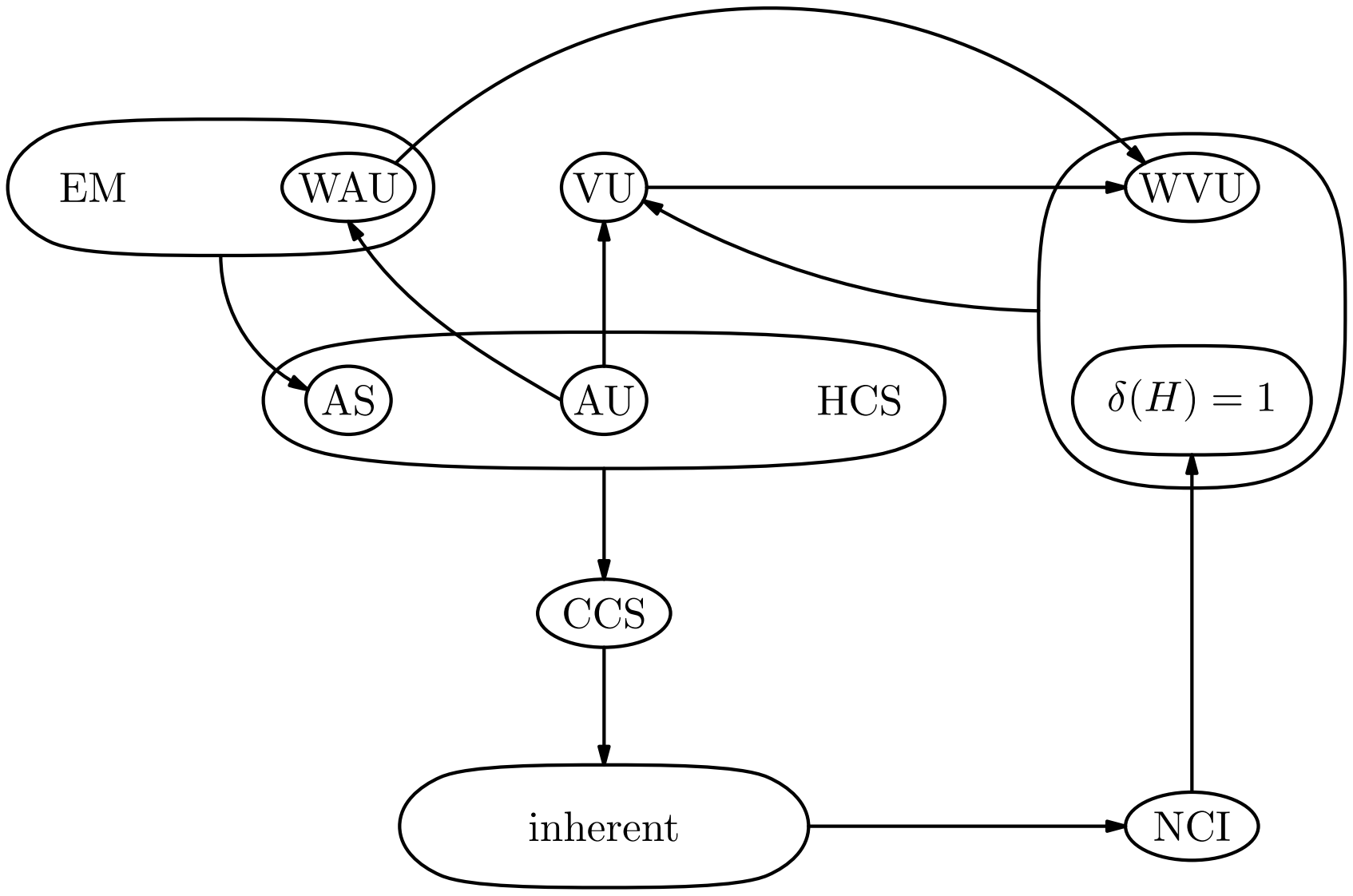

  \centering
  \mpfile{diag}2
  \caption{Implications between various properties of connected two-rooted graphs.
  In general, the implications are of Horn-type, that is, several properties imply another.}
  \label{fig:implications}
\end{figure}

\Cref{fig:implications} summarizes the implications between various properties of two-rooted graphs studied in this paper.
For convenience of the reader, we collect in Appendix \ref{sec:dictionary} abbreviations of various properties used in this paper.

\section{Inherent paths}\label{sec:paths}

Here we extend \cref{thm:endpoints} to the NNI setting and prove the following variant of \cref{thm:paths_scc}, 
thus completely characterizing inherent two-rooted paths.

\begin{thm}\label{thm:ll'}
For every two-rooted path $\cH$, the following properties are equivalent.
\begin{enumerate}[(a)]
    \item\label{item-a} $\cH$ is inherent.
    \item\label{item-b} $\cH$ satisfies the NCI property.
    \item\label{item-c} $\cH$ is $T_1(l,0,l')$ for some $l\ge 1$ and $l'\ge 0$.
    \item\label{item-d} $\cH$ is CCS.
\end{enumerate}
\end{thm}

\begin{rem}\label{rem:terminology}
    Recall that the concepts of inherence differ in the induced~\cite{gurvich2023avoidability} and NNI settings. 
    In particular, it was proved in~\cite{gurvich2023avoidability} that $T_1(2,0,2)$ is confined by the Petersen graph in the induced setting.
\end{rem}

The implication $\eqref{item-a} \Rightarrow \eqref{item-b}$  holds for general connected two-rooted graphs (cf.~\cref{obs:NCI-NNI}), while the implications and $\eqref{item-d} \Rightarrow \eqref{item-a}$ and $\eqref{item-b} \Rightarrow \eqref{item-c}$ are immediate consequences of the definitions.
Hence, it remains to show that every two-rooted graph of the form $T_1(l,0,l')$, with $l\ge 1$ and $l'\ge 0$, is CCS.
In subsequent analysis we denote the number of vertices by $n=l+l'$. 
In the case $n = 1$ the implication follows from a well-known fact that in every connected graph with at least two vertices there exists a vertex such whose removal results in a connected graph. 
So, from now on we assume that $n\ge 2$.
We prove the implication by showing that $T_1(l,0,l')$ is AS and AU, as well as HCS, in the course of \cref{lm:ll'-AS-AU,lm:ll'-HCS}, and then applying \cref{lm:AS-AU-HCS->CCS}.

\begin{lem}\label{lm:ll'-AS-AU}
  $T_1(l,0,l')$ is AS and AU.
\end{lem}

\begin{proof}
  The extension graph of $T_1(l,0,l')$ contains a path $P_{n+1}$. 
  Thus it contains at least two paths $P_n$. 
  So,  $T_1(l,0,l')$ is EM.

  Next, we prove that $T_1(l,0,l')$ is WAU. 
  For this purpose we fix a $T_1(l,0,l')$-smooth graph~$G$ and use the $G$-equivalence relation on arcs of $T_1(l,0,l')$. 
  Let $a_1,\dots,a_{n}$ be the vertices of $T_1(l,0,l')$ and $e_i =\{a_i, a_{i+1}\}$, $1\leq i\leq n-1$ be its edges. 
  Arcs of $T_1(l,0,l')$ are $e^+_i=(a_i,a_{i+1})$ and $e^-_i=(a_{i+1},a_{i})$, $1\leq i\leq n-1$.
  
  Let $P=(x_1,\dots,x_{n})$ be a copy of $P_n$ in a $T_1(l,0,l')$-smooth graph $G=(V,E)$. 
  Note that $P'=(x_n,\dots,x_{1})$ is also a copy of $T_1(l,0,l')$ in $G$. 
  Since $G$ is $T_1(l,0,l')$-smooth, both endpoints $x_1$, $x_{n}$ are roots of non-simplicial copies of $T_1(l,0,l')$, both of them have neighbors in $V\setminus V(P)$, say $x_0x_1\in E$, $x_{n+1}x_n\in E$ such that $x_0\ne x_i$, $x_{n+1}\ne x_i$ for $1\leq i \leq n$ (it is possible that $x_0=x_{n+1}$ if $l'\ne 0$).

  Looking at the path $(x_0,\dots,x_{n-1})$ we see  $e^{+}_i\sim_G e^+_{i+1}$ and $e^{-}_i\sim_G e^-_{i+1}$  for $1\leq i <n-1$. 
  So all $e^+_i$ are $G$-equivalent as well as all $e^-_i$.
  Now look at the copy  $P'=(x_n,\dots,x_{1})$ of $T_1(l,0,l')$ in $G$. 
  Comparing $P$ and $P'$ shows that $e^+_1\sim_G e^{-}_{n-1}$ and $e^{-}_{1}\sim_G e^+_{n-1}$.
Thus every two arcs of $T_1(l,0,l')$ are $G$-equivalent to each other and, due to \cref{lm:equivalence}, $T_1(l,0,l')$ is WAU and thus also WVU.
\cref{lm:EMandWAU->AS} implies that $T_1(l,0,l')$ is also AS. 

To conclude the proof we now show that $T_1(l,0,l')$ is AU. 
Since $n\ge 2$, by \cref{lm:H-cover}, $T_1(l,0,l')$ is VU. 
Let $G=(V,E)$ be an $T_1(l,0,l')$-smooth graph, $(x,y)\in E$ be an arc.
By VU, there exists a copy $(x_1,\dots,x_{n}) $ of $T_1(l,0,l')$ such that $x_1=x$. 
We distinguish the following cases:
\begin{itemize}
    \item If $y=x_2$, we are done.
    \item If $y = x_i$ for some $i>2$, the copy $(x_{i-1}, x_{i-2}, \dots, x_1, x_i, x_{i+1}, \dots ,  x_{n-1}, x_n)$ contains the arc $(x,y)$.
    \item Otherwise, the copy $(x_{n-1},\dots, x_{1},y)$ contains the arc $(x,y)$.   
\end{itemize}
So, in all cases the arc $(x,y)$ is in a copy of $T_1(l,0,l')$. 
Now AU follows from the fact that $T_1(l,0,l')$ is WAU.
\end{proof}

\begin{lem}\label{lm:ll'-HCS}
  $T_1(l,0,l')$ is HCS.
\end{lem}

\begin{proof}
  Let $G=(V,E)$ be an  $T_1(l,0,l')$-smooth graph, $(U,R)$ be an inclusion-maximal \hbox{$P_n$-avoiding pair} in $G$, and $P=(a_1,\dots,a_{n})$ be a copy of $P_n$ in $G[R]$. 
  Since $G$ is $T_1(l,0,l')$-smooth, both endpoints $a_1$, $a_{n}$ of $P$ are roots of non-simplicial copies of $T_1(l,0,l')$, hence, both of them have neighbors outside of $V(P)$. 
  We now show that either $N_G(\{a_1, a_n\})\setminus V(P)\subseteq U $ (in which case $P$ is an \emph{outer} path) or $\{a_1, a_n\}\cap N_G(U)=\es$ (in which case $P$ is an \emph{inner} path).
 
  Suppose that $a_1\notin N_G(U)$ and $a_{n}\in N_G(U)$. 
  An extension of a copy of $T_1(l,0,l')$ gives an edge $a_0a_1\in E$ such that $a_0\in R\setminus V(P)$. 
  We get a copy $(a_0,a_1,\dots,a_{n-1})$ of $T_1(l,0,l')$ in $R_{a_n}$, a contradiction with the hitting condition.
  Similarly, we obtain a contradiction if $a_n\notin N_G(U)$ and $a_{1}\in N_G(U)$.
  Hence, either none of $a_1, a_n$ has a neighbor in $U$ or both of them have neighbors in $U$. 
  In the second case neither of $a_1,a_n$ has a neighbor in $R\setminus V(P) $ due to the same argument.
  Consequently, $P$ is an inner path in the former case and an outer path in the latter case. 
  Note that we have proved that if one of the endpoints of $P$ has a neighbor in $U$, then $P$ is outer, and if not, it is an inner path. 
   
  Assume that $P$ is an outer path, $a_1, a_n\in N_G(U)$. To complete this case, let us prove by contradiction that $R\setminus V(P)=\es$ (which is equivalent to $|R| = n$). Suppose that $R\setminus V(P)\ne\es$.
  Due to the hitting condition, there is no edge in $G$ connecting a vertex in $U$ with a vertex in $R\setminus V(P)$.
  This implies that $V(P)$ is a cut set in $G$.
  Hence, there exists an edge $a_k r$ such that $r\in   R\setminus  V(P)$ and $k\in [n]$. 
  Let $k$ be the smallest number satisfying this condition.   Since $P$ is an outer path, $1< k < n$. 
  Choose  $(v_1,  v_2, \dots, v_q)$ as a longest path in $G[(R\setminus V(P))\cup   \{ a_k\}] $ starting at $ a_k = v_1$.  
  Due to maximality of $U$, $q<k$: otherwise, the path 
  $(v_q, v_{q-1},\dots, v_1,a_{k+1},\dots a_{n})$ contains a subpath with $n$ vertices that does not contain $a_1\in N_G(U)$, which is impossible due to the hitting condition.  
  Since  $ a_1$, $a_{n} $  have neighbors $u', u''$ in $U$ and $G[U]$ is connected, there exists a cycle $( a_1,\dots,  a_{n},  a_{n+1},\dots a_{n+j},  a_1)$, where
  $a_i\in U$ for all $n< i\leq n+j$ and $u'= a_{n+j}$, $u'' =a_{n+1} $. 
  Take a copy $F =  (v_q, v_{q-1},\dots, a_{k}, a_{k+1},\dots,  a_{k+n-q})$ of $P_n$ (indices in $a_i$ are taken modulo  $n+j$). 
  Furthermore, this copy is simplicial in $G$: there are no edges from $v_q$ to
  $U$ due to the hitting condition, there are no edges from $v_q$ to   $ a_i$, $i<k$, by the choice of $k$,
  for all $k\leq i\leq n$ the vertex $a_i$ belongs to $F$  since $q<k$, and, finally, there are no edges from   $v_q$ to $R\setminus V(P)$ since the path $v_1,  v_2, \dots, v_q$  is a longest one. We get a contradiction, since $G$ has no simplicial copies of $T_1(l,0,l')$.

By the previous arguments, we may assume that all the copies of $P_n$ in $G[R]$ are inner paths.
To complete the proof, we verify that if there exists a vertex $x\in N_G(U)$ such that $R_x$ is connected, then $|R|=|V(H)|$.
Fix $x\in N_G(U)$ such that $R_x$ is connected. 
Due to the hitting condition, $x$ is contained in every copy of $P_n$ in $G[R]$.
Choose the minimal $q$ such that $x=x_q$ in a copy $(x_1, \dots, x_n)$ of $P_n$ in $G[R]$. 
Suppose that $q>1$. 
Take an extension $(x_{n+1}, x_n, \dots, x_1; z)$ of a copy $(x_n,\dots, x_1)$, where $x_{n}x_{n+1}$ and $x_{n-l+1}z$ are extension edges. 
If $x_{n+1}\in U$, then the path $(x_n,\dots, x_1)$ is outer, a contradiction.
Otherwise, since $q\ge 2$, the copy $(x_2,\dots, x_{n+1})$ of $P_n$ contains the vertex $x = x_q$, moreover, this vertex corresponds to the $(q-1)$-the vertex of the path, a contradiction with the definition of $q$. 
Thus $q=1$ and $x=x_1$. 
Similar arguments show that, since $(x_1, x_2, \dots, x_n)$ is an inner path, we come to a contradiction with the hitting condition. 
\end{proof}

\cref{thm:paths_scc} follows from \cref{lm:ll'-AS-AU,lm:ll'-HCS} applied to \cref{lm:AS-AU-HCS->CCS}.

\section{Inherent forks} \label{sec:fork}

\begin{sloppypar}
Recall that the \emph{fork} $F(1,1,2)$ (or, equivalently, $F(2,1,1)$) is the unique $5$-vertex tree with maximum degree $3$.
We denote a fork with vertex set $\{v_1,\ldots, v_5\}$ and edge set $\{v_1v_2,v_2v_3,v_3v_4,v_3v_5\}$ by $v_1v_2v_3v_4v_5$.
We prove in \cref{subsec:forks} that there are exactly three two-rooted fork graphs that are PE-inherent.  
They are shown in  \cref{pic:T1-211,pic:T3-11,pic:T1-112}.
For later use, we identify the vertices of the graphs  as shown in this figure and denote a copy of the corresponding graph by the sequence shown in the figure.
Below we refer to these two-rooted graphs as \emph{PE-forks}.
\end{sloppypar}

 \begin{figure}[!h]
   \centering
   \begin{minipage}{0.3\textwidth}
   \centering
  \mpfile{fork}{1}
  \caption{$(s,t,a,b,c)$ is a copy of the two-rooted fork  $T_1(2,1,1)$.}\label{pic:T1-211}
   \end{minipage}
   \quad
  \begin{minipage}{0.3\textwidth}
   \centering
  \mpfile{fork}{2}
  \caption{$(c,t,a,b,s)$ is a copy of the  two-rooted fork $T_3(1,1)$.}\label{pic:T3-11}
  \end{minipage}
\quad
\begin{minipage}{0.3\textwidth}
   \centering
  \mpfile{fork}{3}
  \caption{$(d,c,a,b,s)$ is a copy of  the two-rooted fork $T_1(1,1,2)$.}\label{pic:T1-112}
  \end{minipage}
  \end{figure}

\begin{thm}\label{thm:forks}
  PE-forks are CCS graphs. 
  Hence, they are inherent.
\end{thm}

Due to \cref{lm:AS-AU-HCS->CCS}, to prove \cref{thm:forks} it is enough to prove that a PE-fork is AS, AU, and HCS.  

\begin{lem}\label{lm:fork-AS-AU}
  A PE-fork $\cF$ is AS and AU.
\end{lem}

\begin{proof}
In \cref{pic:fork-WAU}(b-d) the extension graphs for PE-forks are shown, together with the labels of vertices that we use in this proof.

It is easy to verify  that PE-forks are EM: for $T_1(2,1,1)$ an additional copy is, e.g., $(b,a, t,s,v_2)$;  for $T_3(1,1)$ an additional copy is, e.g., $(s,a, t,c,v_2)$; for $T_1(1,1,2)$ an additional copy is, e.g., $(v_2,s,a, b,c)$.

\begin{figure}[!h]
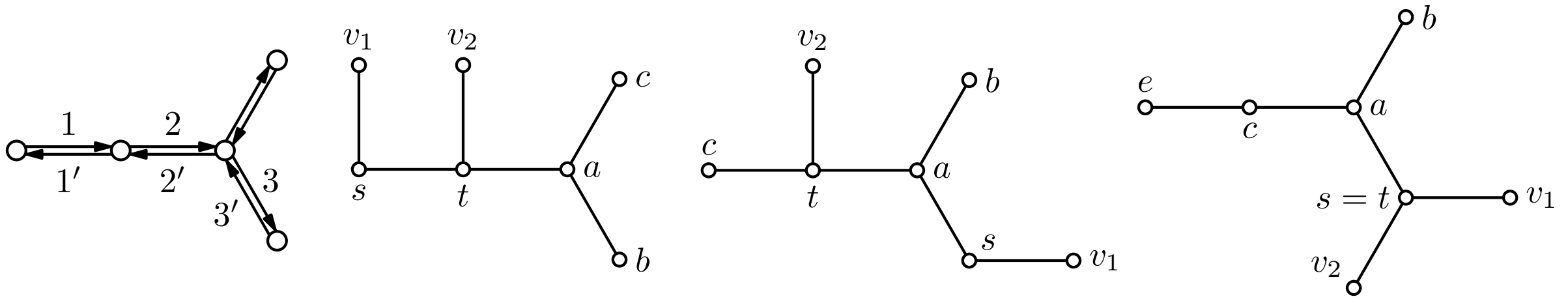

    \centering
    \begin{tabular}{c@{\quad}c@{\quad}c@{\hskip4pt}c}
   \begin{minipage}{0.18\textwidth}
   \centering
   \mpfile{fork}{10}\vspace*{7pt}
   \end{minipage}
   &
   \begin{minipage}{0.2\textwidth}
   \centering
   \mpfile{fork}{11}\vspace*{7pt}
   \end{minipage}
   &
   \begin{minipage}{0.26\textwidth}
   \centering
   \mpfile{fork}{21}\vspace*{7pt}
   \end{minipage}
   &
   \begin{minipage}{0.26\textwidth}
   \centering
   \mpfile{fork}{31}\vspace*{7pt}
   \end{minipage}
   \\ 
   (a)&(b) $T_1(2,1,1)$& (c) $T_3(1,1)$& (d) $T_1(1,1,2)$
    \end{tabular}
   \caption{Naming arcs of the fork; and the extension graphs of PE-forks.}\label{pic:fork-WAU}
 \end{figure}

Next, we prove that $\cF$ is WAU. For this purpose we use \cref{lm:equivalence}. 
So, we are going to prove that for every $\cF$-smooth graph $G$, all arcs of $\cF$ are $G$-equivalent.  
Name (some of the) arcs of the fork as shown in \cref{pic:fork-WAU}(a).  
Note that it is enough to consider only three pairs of arcs due to automorphisms of the fork, since arcs in the same orbit of the automorphism group are obviously equivalent. 
    
Now we compare arcs in the original copy and additional copies in its extensions in $G$.  

    Let $\cF= T_1(2,1,1)$ and $G$ be a $T_1(2,1,1)$-smooth graph.
    The~copy $(b,a, t, s, v_2)$ of $T_1(2,1,1)$ shows that $1\sim_G 3'$, $1'\sim_G3$,
    $2\sim_G2'$.
    Furthermore, the~copy $(v_1,s,t,a,v_2)$ gives equivalences $1\sim_G2$,  $1'\sim_G2'$, $2\sim_G3$, $2'\sim_G3'$.
    Hence, all arcs of $\cF$ are $G$-equivalent: \[3\sim_G 1'\sim_G2'\sim_G 2\sim_G 1\sim_G 3'\,.\]
We conclude from \cref{lm:equivalence} that  $T_1(2,1,1)$ is WAU.

    Let $\cF=T_3(1,1) $ and $G$ be a $T_3(1,1)$-smooth graph.
    The~copy $(b,a, t,c, v_2)$ of $T_3(1,1)$ shows that  $1\sim_G3'$, $1'\sim_G3$,  $2\sim_G2'$, see \cref{pic:fork-WAU}(c). 
    The copy $(v_1,s, a, t,b)$  shows that $3\sim_G2'$, $3'\sim_G2$. So, all arcs of $T_3(1,1)$ are $G$-equivalent:
\[
1'\sim_G 3\sim_G 2'\sim_G 2\sim_G 3'\sim_G 1.
\]
From \cref{lm:equivalence} we conclude that $T_3(1,1)$ is WAU.

Finally, we consider the case  $\cF = T_1(1,1,2)$. This case is more complicated.  
Let $G$ be a $T_1(1,1,2)$-smooth graph.
The copy $(v_2,s,a,b,c)$ from the extension graph in \cref{pic:fork-WAU}(d) gives equivalences $3\sim_G2'$, $3'\sim_G2$. 
One more copy of $T_1(1,1,2)$ in the extension graph gives the same equivalences.  
So, we take an extension of $(v_2,s,a,b,c)$, which gives an edge $cv_3$, $v_3\ne e$ and  a new copy $(s,a,c,e,v_3)$. 
The copy gives equivalences $3\sim_G1'$, $3'\sim_G1$, $2\sim_G2'$. 
So, all arcs of $T_1(1,1,2)$ are $G$-equivalent:
\[
1'\sim_G 3\sim_G 2'\sim_G 2\sim_G 3'\sim_G 1,
\]
and an application of  \cref{lm:equivalence} shows that  $T_1(1,1,2)$ is WAU.

Since PE-forks are EM and WAU , they are AS due to \cref{lm:EMandWAU->AS}.

It remains to prove that PE-forks are AU.   
Let $\cF$ be a PE-fork and $G=(V,E)$ be an $\cF$-smooth graph. 
Take an arc $(x,y)$ in $G$. 
Since $\cF$ is WAU, it is WVU by a direct implication and VU by \cref{lm:H-cover}. 
Take a copy of $\cF$ containing $x$. 
If the copy contains also $xy$, we are done, since $\cF$ is WAU. 
If not, there should be an edge $xz$ in the copy. 
Since $\cF$ is WAU, this implies that there exists a copy of $\cF$ such that the arc $(x,z)$ corresponds to the arc from the vertex of degree $2$ to the vertex of degree $3$ in the fork. 
If this copy contains the edge $xy$, we are done. 
If the copy does not contain $y$, then substitute a pendant vertex adjacent to $x$ by $y$, and we are also done. 
The remaining case is unique up to isomorphism, it is shown in \cref{pic:fork-AU}.
  \begin{figure}[!h]
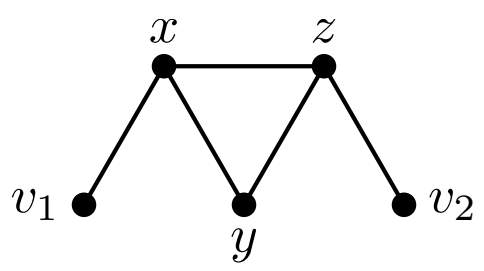

    \centering
    \mpfile{fork}{4}
    \caption{The special case in the proof of AU for a PE-fork.}\label{pic:fork-AU}   
  \end{figure}
In this case $xy$ is an edge of the copy $(v_2,z,x,y,v_1)$. Since $\cF$ is WAU, this completes the proof that $\cF$ is AU.
\end{proof}

To prove that PE-forks are HCS, we will need the following
claim about \emph{fork-free graphs}, that is, graphs containing no copies of the fork.

\begin{lem}\label{fork-free}
  Any connected fork-free graph on at least $5$ vertices is either a~star $K_{1,n}$, a path $P_n$, or a cycle $C_n$ for some $n$.
\end{lem}

\begin{proof}
Let $G = (V,E)$ be a connected fork-free graph with $|V|\ge 5$.
  If all vertices in $G$ have degree at most $2$ then $G$ is path or a cycle, since $G$ is connected. 
  From now on assume that there exists $v\in V$ such that $d_G(v)\geq3$. 
  We need to prove that $G \cong K_{1,n}$ (for $n = |V|-1$).
  
  Suppose the contrary. 
  Since $G$ is connected, there exists an edge $xy\in E$, $x\in N(v)$, $y\ne v$.  
  If $d_G(v)=3$, then, since $G$ has at least $5$ vertices, we may assume that $y$ is chosen outside $N(v)$, and we obtain a copy $yxvv_1v_2$ of the fork, where $v_1,v_2\in N(v)$ with $v_1\neq v_2$. 
  If $d_G(v)\geq 4$, then choose $v_1,v_2\in N(v)\setminus \{x,y\}$ and get a copy $yxvv_1v_2$ of the fork. 
  Thus, $G$ is not fork-free, a contradiction.
\end{proof}

\begin{lem}\label{lm:fork->HCS}
  A PE-fork $\cF$ is HCS.
\end{lem}

\begin{proof}
Let $G=(V,E)$ be an $\cF$-smooth graph, and $(U,R)$ be an inclusion-maximal $\cF$-avoid\-ing pair $(U,R)$ in $G$, and   $x\in N_G(U)$ such that $R_x$ is connected.  
Since $\cF$ is AU by \Cref{lm:fork-AS-AU}, $\cF$ is also VU.
Hence, by \cref{cor:deg3}, we have $\delta(G)\geq3$. 

We are going to prove that $|R| = |V(\cF)|=5$. 
So we need to exclude the case when $|R|> 5$. 
Note that in this case $|V(R_x)|\geq 5$ and $R_x$ is fork-free due to the hitting condition.
From \cref {fork-free} we conclude that there are exactly three possible cases for $R_x$.

\textbf{\boldmath Case $R_x= K_{1,n}$, $n\geq4$.} 
Let $v_i$ be pendant vertices of the star $R_x$.
Then $d_G(v_i)\geq \delta(G) \geq3$. 
This implies that  $v_i\in N_G(U)$ for all $i$, since there are only two possible edges in $G[R]$ incident to $v_i$: the edge $v_iz$, where $z$ is the center of the star and (possibly) the edge $v_ix$. 
Therefore $\{x, v_1,v_2,v_3,v_4\}\subseteq N_G(U)$.
By the hitting condition, each copy of of the fork in $G[R]$  contains $N_G(U)$, and therefore $G[\{x, v_1,v_2,v_3,v_4\}]$ is a fork.
We come to a contradiction, since $G[\{x, v_1,v_2,v_3,v_4\}]$ has an independent set of size 4 and the fork does not.

\textbf{\boldmath Case $R_x= P_{n}$, $n\geq5$.} 
Enumerate vertices along the path: $v_1, v_2,\dots, v_{n}$, where $n\geq5$. 
Since $d_G(v_i)\geq \delta(G) \geq3$, $\{v_1,v_{n}\}\in N_G(U)$. Furthermore, since any copy of the fork should contain the set $N_G(U)$, the graph $R_{v_1}$ is fork-free and, since the fork has only $5$ vertices, there exists a vertex $v_i$ in the path such that $v_i\notin N_G(U)$ and therefore  $xv_i\in E$. 
Since $n\geq5$ and the path is induced, $v_i$ is not adjacent to at least one endpoint of the path. 
By symmetry, we may assume that $i>2$.
Notice that the connected component $R'$ of the graph $R_{v_1}$ that contains the vertex $v_i$ satisfies $|V(R')|\geq 5$, since $R'$ contains all the vertices $v_j$ for $j>1$, as well as vertex $x$. 
Since $v_i$ has degree $3$ in $R'$, we conclude that $R'$ is neither a star, nor a path, nor a cycle. 
Since $R'$ is fork-free, this contradicts \cref{fork-free}.

\textbf{\boldmath Case  $R_x=C_n$, $n\geq5$.} 
Again, there exists a vertex $v$ in the cycle such that $v\notin N_G(U)$. 
Since $d_G(v)\geq \delta(G) \geq3$, we conclude that $vx\in E$. 
Enumerate vertices along the cycle starting from $v$: $v_1, v_2,\dots, v_{n}$, $v_1=v$, where $n\geq5$. 
For $3\leq i\leq n-1$, in the connected component $R'$ of the graph $R_{v_i}$ that contains the vertex $v_1$, the vertex $v_1$ has degree at least~$3$ (due to the edges $v_1x$, $v_1v_2$, and $v_1v_{n}$) and, moreover, $|V(R')|\geq 5$, since $R'$ contains at least four vertices of the cycle, as well as vertex $x$. 
Since $v_1$ has degree at least $3$ in $R'$, we conclude that $R'$ is neither a star, nor a path, nor a cycle. 
By \cref{fork-free}, $R'$ contains a copy of the fork and, hence, so does $R_{v_i}$.
By the hitting condition, this implies that $v_i\notin N_G(U)$, hence, the only possible neighbors of $v_i$ in $G$ are $x$, $v_{i-1}$, and $v_{i+1}$. 
Since also $d_G(v_i)\geq \delta(G) \geq3$, we conclude that $d_G(v_i)=3$ and, moreover, $xv_i\in E$  for $3\leq i\leq n-1$. 
This leads to a contradiction with the fact that $G$ is $\cF$-smooth, since a simplicial copy of $\cF$ does exist under these conditions, see \cref{pic:Cn}, where simplicial copies are shadowed.
  \begin{figure}[h]
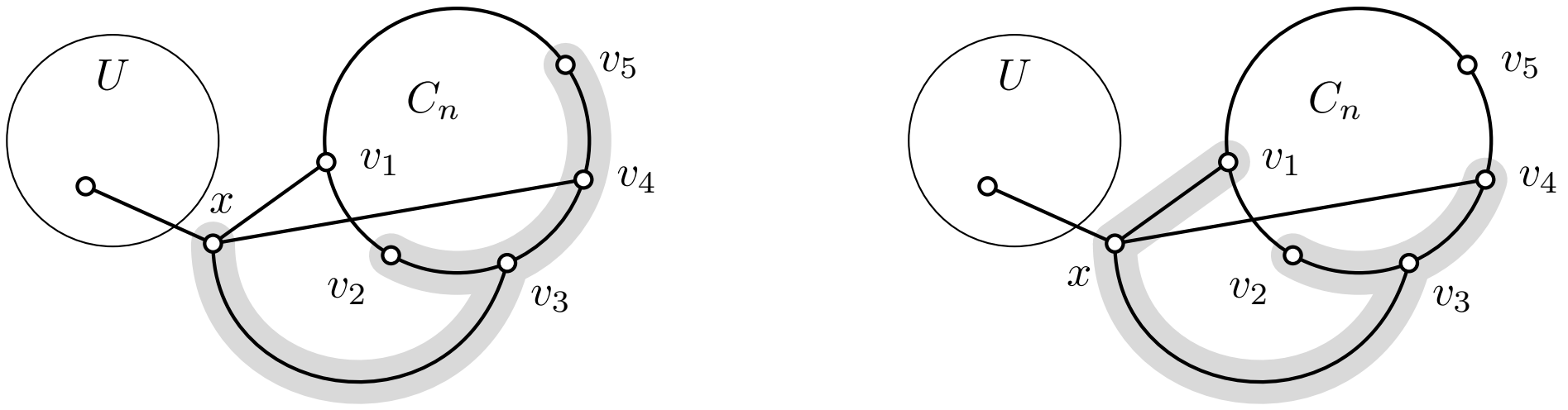

    \centering
    \mpfile{forks}{60}\hspace*{2cm}\mpfile{forks}{61}
    \caption{Simplicial copies: $(v_5,v_4,v_3,v_2,x)$ for
      $T_1(2,1,1)$ and for $T_3(1,1)$ (on the left); $(v_1,x,v_3,v_2,v_4)$ for
      $T_1(1,1,2)$   (on the right).}\label{pic:Cn}
  \end{figure}
\end{proof}

\cref{thm:forks} follows from \cref{lm:fork-AS-AU,lm:fork->HCS} applied to \cref{lm:AS-AU-HCS->CCS}.

\section{Proof of \cref{thm:T1212}}\label{sec:T1212}

Let recall the statement of the theorem.
\thmsixvertex*

The~proof of \cref{thm:T1212} is similar to the proofs for PE-forks. 
In fact, we conjecture that $T_1(2,1,2)$ is CCS; however, the analysis is much more involved.
Inherence is easier to prove than CCS, since we can stop our analysis whenever an avoidable copy is detected.

\begin{figure}[!h]
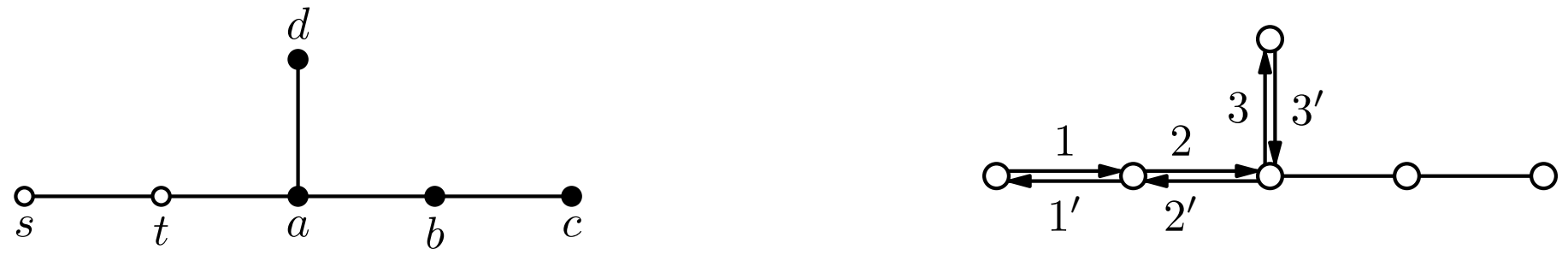

    \centering
    \begin{tabular}{cc}
   \begin{minipage}{0.43\textwidth}
   \centering
   \mpfile{T1212}{4}\vspace*{7pt}
   \end{minipage}
   &
    \begin{minipage}{0.43\textwidth}
   \centering
   \mpfile{T1212}{5}\vspace*{7pt}
   \end{minipage}
  \\ 
  (a) $(s,t,a,b,c,d)$ is a copy of   $T_1(2,1,2)$
  &(b) Names of the  arcs of  $T_1(2,1,2)$
    \end{tabular}
   \caption{Naming arcs and vertices of $T_1(2,1,2)$.}\label{pic:T1212-names}
 \end{figure}

We will indicate  a copy of $T_1(2,1,2)$ by a sequence of vertices in the order shown in \cref{pic:T1212-names}(a). 
For the purpose of analyzing $G$-equivalence with respect to some \hbox{$T_1(2,1,2)$-smooth} graph $G$, we name arcs of $T_1(2,1,2)$ as shown in \cref{pic:T1212-names}(b). 
Due to automorphisms of $T_1(2,1,2)$, it is enough to consider only three pairs of arcs.

 \begin{lem}\label{T1212isWAU}
   $T_1(2,1,2)$ is WAU.
 \end{lem}
 \begin{proof}
Let $G$ be a $T_1(2,1,2)$-smooth graph. 
We fix a copy $(s,t,a,b,c,d)$ of $T_1(2,1,2)$ in $G$ and an extension of it using edges $sv_1$ and $tv_2$, as shown in \cref{pic:T1212-WAU}.

 \begin{figure}[!h]
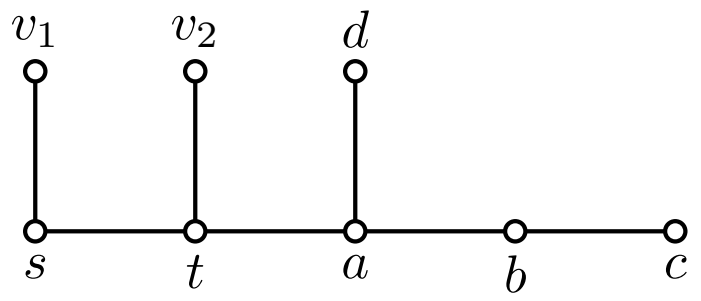

    \centering
    \mpfile{T1212}{501}
    \caption{An extension of $T_1(2,1,2)$.}\label{pic:T1212-WAU}
 \end{figure}

Let us compare the initial copy of $T_1(2,1,2)$ and the copy $(v_1,s,t,a,d, v_2)$ in the extension graph.
This gives equivalences $1\sim_G2$, $1'\sim_G2'$, $2\sim_G2'$, $3\sim_G1'$, $3'\sim_G1$. 
Thus, all arcs of $T_1(2,1,2)$ are $G$-equivalent:
\[
3\sim_G1'\sim_G2'\sim_G 2\sim_G1\sim_G 3',
\]
and an application of \cref{lm:equivalence} shows that $T_1(2,1,2)$ is WAU.
\end{proof}

 \begin{cor}
    Let $G$ be a $T_1(2,1,2)$-smooth graph. 
    Then $\delta(G)\geq3$.
 \end{cor}
 
 \begin{proof}
  Recall that WAU implies WVU. 
  The corollary follows from \cref{T1212isWAU,lm:H-cover,cor:deg3}.
 \end{proof}

 \begin{lem}\label{U>=2}
   Let $G$ be a $T_1(2,1,2)$-smooth graph. Then there exists an   inclusion-maximal  avoiding pair $(U,R)$ for  $T_1(2,1,2)$ in $G$ such that $|U|\geq2$.
 \end{lem}
 
 \begin{sloppypar}
 \begin{proof}
   Take a copy of $T_1(2,1,2)$ and its extension graph shown in \cref{pic:T1212-WAU}. 
   The copy $(v_2,t,a,b,c,d)$ certifies that $(\{v_1,s\}, V(G)\setminus(\{v_1,s\}) $ is an avoiding pair for  $T_1(2,1,2)$ in $G$. 
   It can be extended to an inclusion-maximal avoiding pair $(U,R)$ and the size of $U$ is at least $2$ in that pair.
 \end{proof}
\end{sloppypar}

\begin{proof}[Proof of \cref{thm:T1212}.]
Suppose for a contradiction that $T_1(2,1,2)$ is not inherent.
Then, there exists a connected confining graph $G$ for $T_1(2,1,2)$.
In particular, the graph $G$ is \hbox{$T_1(2,1,2)$-smooth}.
   By \Cref{U>=2}, there exists an inclusion-maximal avoiding pair $(U,R)$ for $T_1(2,1,2)$ in $G$ such that $|U|\geq2$. 
   Due to maximality of $U$, the hitting condition holds: for each $x\in R\cap N_G(U)$, the graph $R_x= G[R\setminus\{x\}]$ is  $T_1(2,1,2)$-free. 

   Let $A = (x_1, y_1, z_1, y_2, x_2, z_2)$ (see \cref{pic:T1212-names}(a) for naming convention) be a copy of $T_1(2,1,2)$ contained in $G[R]$.
    Since $G$ is a confining graph for $T_1(2,1,2)$, this copy of $T_1(2,1,2)$ is not avoidable in $G$; hence, there exists a~non-closable extension by edges $x_1w_1$, $y_1w_2$. 
   Let $R'=R\setminus V(A)$. 
   Recall that $G[U]$ is connected. 
   So, if $w_1,w_2\in U$, then the extension is closable. 
   Note that the case $w_1\in U$, $w_2\in R'$ is impossible because in this case $R_{x_1}$ does contain a copy $(w_2, y_1, z_1, y_2, x_2, z_2)$ of  $T_1(2,1,2)$, in contradiction with the hitting condition. 
   To complete the proof, we need to consider two remaining cases.
 
   \textbf{\boldmath Case $w_1\in R'$, $w_2\in R'$.} 
   First, we note that $x_1\notin N_G(U)$, since otherwise we come to a contradiction with the hitting condition as explained above. 
   From the hitting condition we also conclude that $x_2,y_2\notin N_G(U)$ since the graphs $R_{x_2}$ and $R_{y_2}$ contain the copy $(w_1,x_1,y_1, z_1, z_2, w_2)$ of $T_1(2,1,2)$.  
   Let $w_1v_1$, $x_1v_2$ be the edges of an extension of this copy.
   Note that it is possible that $\{v_1,v_2\}\cap \{x_2,y_2\}\ne \es$, and  $v_1,v_2\notin U$ since ${(R'\cup\{x_1\})\cap N_G(U)=\es}$. 
   The copy $(v_1,w_1,x_1, y_1, w_2, v_2)$ of $T_1(2,1,2)$ is contained in the graphs $R_{z_1}$ and $R_{z_2}$, see \cref{pic:T1212-Ia} (dashed circles in the figure indicate vertices that are in $R\setminus N_G(U)$). 
   This implies that  $z_1,z_2\notin N_G(U)$. 
   Thus $N_G(U)=\{y_1\}$ since this set is not empty (as $G$ is connected). 

   \begin{figure}[!h]
   \centering
   \begin{minipage}[t]{0.48\textwidth}
   \centering
  \mpfile{T1212}{511}
  \caption{To  the proof that  $ N_G(U)=\{y_1\}$.}\label{pic:T1212-Ia}
   \end{minipage}
   \quad
   \begin{minipage}[t]{0.48\textwidth}
   \centering
  \mpfile{T1212}{512}
  \caption{A contradiction with the hitting condition.}\label{pic:T1212-Ib}
   \end{minipage}
\end{figure}

The copy $(x_2,y_2,z_1,y_1,x_1,z_2)$ of $T_1(2,1,2)$ is not avoidable. 
Due to the previous analysis and the fact that $y_2\notin N_G(U)$, there exists an extension of the copy by edges $x_2v_3$, $y_2v_4$, where $v_3,v_4\in R'$ (see \cref{pic:T1212-Ib}). 
The copy $(z_2,z_1,y_2,x_2,v_3,v_4)$ is contained in $R_{y_1}$ and we come to a contradiction with the hitting condition. 

   \textbf{\boldmath Case $w_1\in R'$, $w_2\in U$.}
   The previous case shows that $w_1\in R'$, $w_2\in U$ holds for every copy of $T_1(2,1,2)$ in $G[R]$.
   Take a copy $A = (x_1, y_1, z_1, y_2, x_2, z_2)$ and the `reflected' copy  $A' = (x_2, y_2, z_1, y_1, x_1, z_2)$. 
From extensions of these copies we conclude that $y_1, y_2\in N_G(U)$.
   
   Next step is to prove that $x_1\in  N_G(U)$. 
   Suppose for contradiction that $x_1\notin  N_G(U)$. 
   Then an extension of the copy $(w_1, x_1, y_1, z_1,z_2, w_2)$ uses edges $w_1v_1$, $x_1v_2$ such that $v_1, v_2\notin U$ (see \cref{pic:T1212-IIa}(a)). 
   Thus the copy $(v_1,w_1,x_1,y_1,z_1,v_2)$ should not be contained in $R_{y_2}$, implying that $y_2\in \{v_1,v_2\}$.
   If $y_2=v_1$, then the copy $(z_2,z_1,y_2,w_1,x_1,x_2)$ is contained in $R_{y_1}$, in contradiction with the hitting condition (see \cref{pic:T1212-IIa}(b)).

   \begin{figure}[!h]
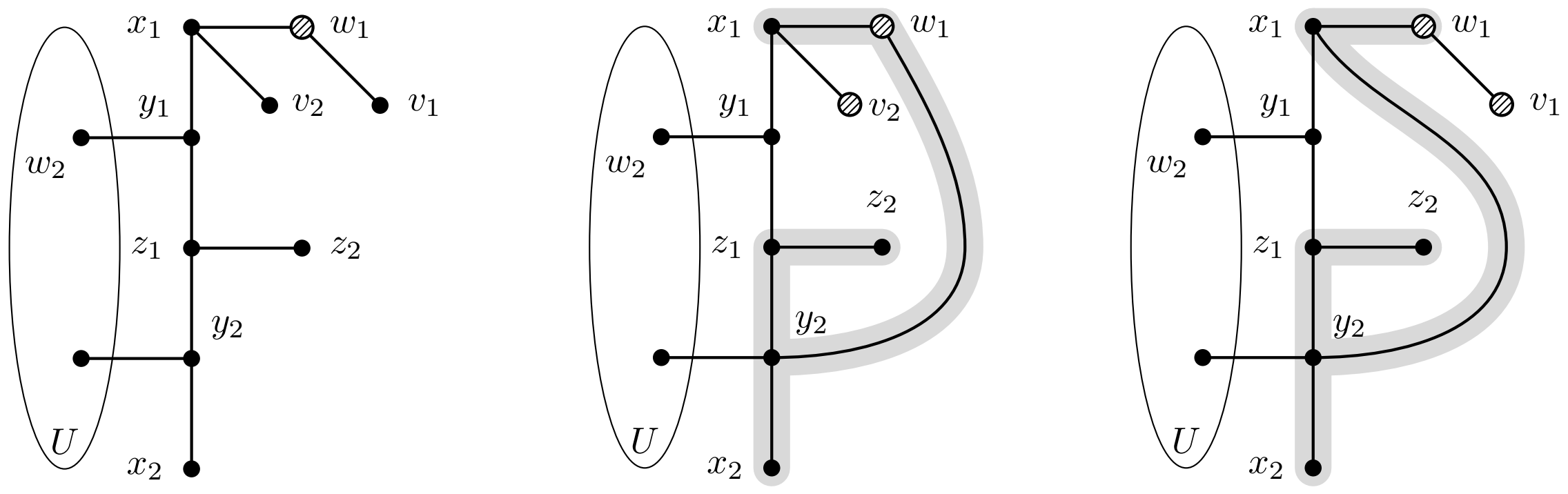

    \centering
    \begin{tabular}{ccc}
   \begin{minipage}{0.3\textwidth}
   \centering
   \mpfile{T1212}{521}\vspace*{7pt}
   \end{minipage}
   &
    \begin{minipage}{0.3\textwidth}
   \centering
   \mpfile{T1212}{522}\vspace*{7pt}
    \end{minipage}
&
    \begin{minipage}{0.3\textwidth}
   \centering
   \mpfile{T1212}{523}\vspace*{7pt}
    \end{minipage}
    \\[10pt]
    (a)& (b) $y_2=v_1$& (c) $y_2=v_2$
    \end{tabular}
    \caption{To the proof that $x_1\in N_G(U)$.}\label{pic:T1212-IIa}
\end{figure}

   If $y_2=v_2$, then the copy $(z_2, z_1, y_2, x_1,w_1,x_2)$ is contained in  $R_{y_1}$ (see \cref{pic:T1212-IIa}(c)).
   Again, we come to a contradiction.
   This shows that $x_1\in  N_G(U)$, as claimed.

Substituting $A$ by $A'$ in the previous argument shows that $x_2\in N_G(U)$ as well.

Neighbors of $x_1, x_2, y_1, y_2$ in $U$ can coincide. 
In any case, since $G[U]$ is connected and $|U|\geq2$, there exists a path $x_1uu'$, where $u,u'\in U$ (see  \cref{pic:T1212-IIb}(a)). 
Thus, the arc $(x_1,w_1)$ is contained in the copy $(u', u, x_1, y_1, z_1, w_1)$ of $T_1(2,1,2) $. 

   \begin{figure}[!h]
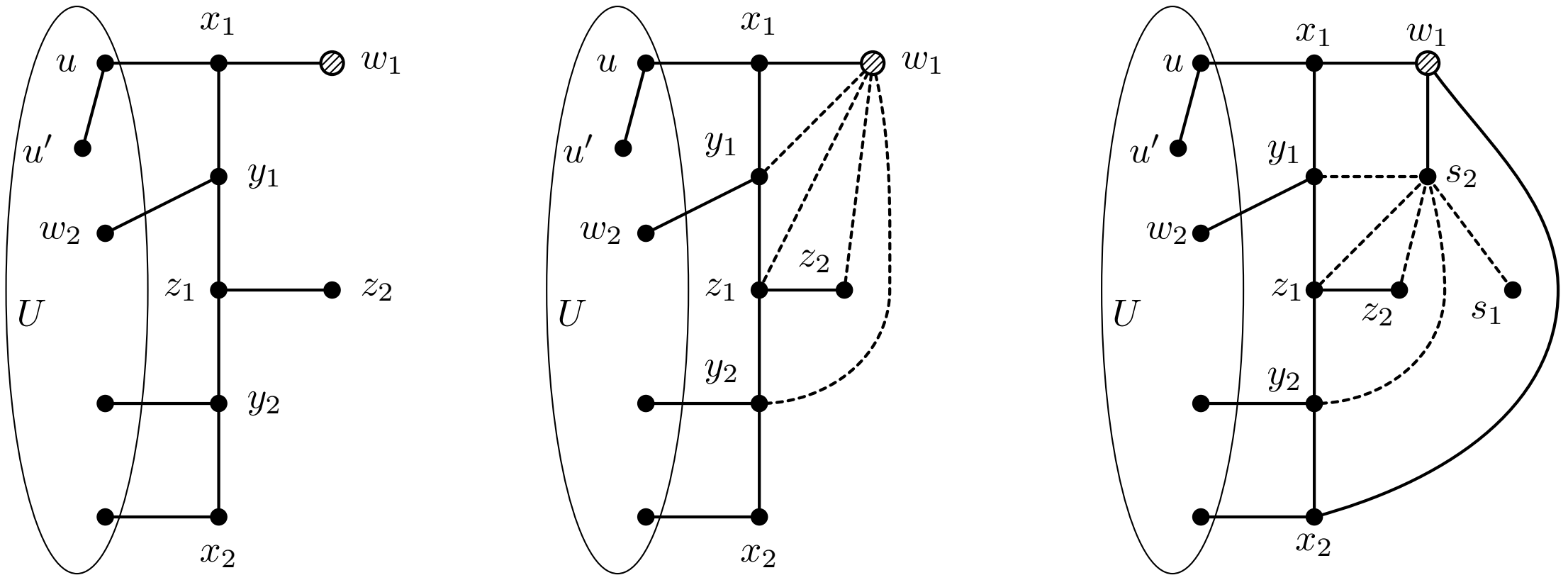

    \centering
    \begin{tabular}{ccc}
   \begin{minipage}{0.28\textwidth}
   \centering
   \mpfile{T1212}{531}\vspace*{7pt}
   \end{minipage}
   &
    \begin{minipage}{0.28\textwidth}
   \centering
   \mpfile{T1212}{532}\vspace*{7pt}
    \end{minipage}
&
    \begin{minipage}{0.33\textwidth}
   \centering
   \mpfile{T1212}{533}\vspace*{7pt}
    \end{minipage}
    \\[10pt]
    (a)& (b) & (c) $s_3=x_2$, $s_2\in R'$
    \end{tabular}
    \caption{Using WAU to complete the  proof. Dashed lines correspond to non-edges.}\label{pic:T1212-IIb}
\end{figure}

Since $T_1(2,1,2) $ is WAU due to \cref{T1212isWAU}, this arc should be an arc corresponding to the arc $(d,a)$ (see \cref{pic:T1212-names}(a)) in a copy $B=(s_1,s_2, w_1, s_3,s_4,x_1)$ of $T_1(2,1,2) $ in $G$. 
Since $w_1\in R'$, neither of $s_2$ and $s_3$ belongs to $U$. 
If both $s_2,s_3\notin V(A)$, then $B$ is contained in $G[R]$. 
Due to the hitting condition, we infer that $\{x_2, y_1, y_2\}\subseteq \{s_1,s_4\}$, which is impossible. 
Thus, at least one of the vertices $s_2,s_3$ should be in $V(A)$.

By symmetry, we may assume that $s_3\in V(A)$. 
A direct case analysis shows that $s_3=x_2$, since in all other cases  we are able to detect a~vertex  $w\in N_G(U)$ and a copy of $T_1(2,1,2) $ contained in $R_w$, in contradiction with the hitting condition, see \cref{pic:T1212-IIb}(b):
\begin{description}
\item [\boldmath $s_3 = y_1$:] $R_{x_1}$ contains  the copy $(w_1,y_1,z_1,y_2,x_2,z_2)$;
\item [\boldmath $s_3 = z_1$:] $R_{y_1}$ contains  the copy $(x_1, w_1,z_1,y_2,x_2,z_2)$;
\item [\boldmath $s_3 = z_2$:] $R_{x_1}$ contains  the copy $(w_1,z_2,z_1,y_2,x_2,y_1)$;
\item [\boldmath $s_3 = y_2$:] $R_{y_1}$ contains  the copy $(z_2,z_1,y_2,w_1,x_1,x_2)$.
\end{description}

This implies that $s_2\in R'$, since otherwise $s_2\in V(A)$, and a symmetric analysis as the one above shows that $s_2=x_2$, implying $s_2 = s_3$, which is impossible.
For all possible variants of placement the vertex $s_1$ we now come to a~contradiction with the hitting condition. 
The case analysis is similar to the previous one, see \cref{pic:T1212-IIb}(c): 
\begin{description}
\item [\boldmath $s_1\in R'\cup\{z_1,z_2\}$:] $R_{y_1}$ contains  the copy $(s_1,s_2,w_1,x_2,y_2,x_1)$;
\item [\boldmath $s_1 = y_1$:] $R_{y_2}$ contains  the copy $(w_1,s_2,y_1,z_1,z_2,x_1)$;
\item [\boldmath $s_1 = y_2$:] $R_{x_1}$ contains  the copy $(w_1,s_2,y_2,z_1, z_2,x_2)$.
\end{description}

We have analyzed all possible cases and come to a contradiction in each of them. 
This completes the proof.   
 \end{proof}

\section{Disconnected inherent graphs}\label{sec:disconn}

Let $H+H'$ denote the disjoint union of graphs $H$ and $H'$.
It seems to us that the following natural statement is true.

\begin{conjecture}
\label{conj:disconnected}
    Let $(H,s,t)$ be a two-rooted graph, and let $H'$ be any graph. 
    Then, the two-rooted graph $(H+H',s,t)$ is inherent if and only if so is $(H,s,t)$.
\end{conjecture}

The backward direction of the \cref{conj:disconnected} holds due to \cite{gurvich2023avoidability}, where the authors proved the statement for the induced case.
The same proof, with some obvious necessary modifications, remains valid for the NNI setting.
In this case, the following analogue of \cite[Proposition 5.8]{gurvich2023avoidability} holds.
We include the proof for completeness.

\begin{prop}\label{thm:adding a component}
Let $(H,s,t)$ be an inherent two-rooted graph and let $H'$ be any graph.
Then the two-rooted graph $(H+H',s,t)$ is inherent.
\end{prop}

\begin{proof}
Let $G$ be a graph containing $H+H'$ as a subgraph.
Fix a copy $\cpi{H + H'}$ of $H+H'$ in $G$,  let $\cpi H$ and $\cpi {H'}$ be the corresponding copies of $H$ and $H'$, and let $X$ be the set of vertices of $\cpi {H'}$.
Observe that $\cpi H$ is a subgraph of $G-X$ isomorphic to $H$.
Thus, since $(H,s,t)$ is inherent, $G-X$ contains an avoidable copy $\cp{H_1}{s_1}{t_1}$ of $(H,s,t)$. 
Then, it can be verified that $(\cpi{H_1}\cup \cpi{H'},\cpi{s_1},\cpi{t_1})$ is an avoidable copy of $(H+H',s,t)$ in $G$.
Thus, $(H+H',s,t)$ is inherent.
\end{proof}

There is a natural way of how one could attempt to prove the forward direction of \cref{conj:disconnected}. 
Suppose that $s,t\in V(H)$ and $(H+H',s,t)$ is inherent.
To show that $(H,s,t)$ is inherent, consider an arbitrary graph $G$ containing
$H$ as a subgraph. Let $G'$ be the disjoint union of $G$ and $H'$. 
Then, $G'$ contains a copy of $H+H'$. 
Since $(H+H',s,t)$ is inherent, $G'$ also contains an avoidable copy $(\cpi{H+H'},\cpi{s},\cpi{t})$ of $(H+H',s,t)$ in $G'$. 
This, however, does not necessarily imply that $\cp Hst$ is an avoidable copy of $(H,s,t)$ in $G$, as illustrated by the following example.
Suppose that $H\cong H'\cong K_3$, with $s$ and $t$ adjacent, and let $G$ be the graph obtained from $K_3$ by adding a pendant edge to each vertex.
Then, there is a simplicial copy of $H+H'$ in $G+K_3$, even though $G$ contains a copy of $H$ but no avoidable copy of $(H,s,t)$.

We point out that there are inherent two-rooted graphs $(H,s,t)$ such that $s$ and $t$ are in different connected components. The following proposition provides an example.

\begin{prop}\label{prop:isolated-root}
Let $H=P_k$, $k>1$, and $(H,s,t)$ be the endpoint-rooted path, which is inherent. 
Let  $H'$ be the graph obtained from $H$ by removing the edge incident with~$s$ (thus, leaving $s$ isolated). 
Then $(H',s,t)$ is inherent.
\end{prop}

\begin{proof}
Let $G$ be a graph containing $H'$  as a subgraph. 
Suppose first that $G$ contains $P_k$ as a subgraph. 
Since $(H,s,t)$ is inherent, $G$  contains an avoidable copy $\cp Hst$ of $(H,s,t)$. 
The copy of $(H',s,t)$ obtained by deleting from $\cp Hst$ the unique edge incident with $\cpi s$ is an avoidable copy of $(H',s,t)$ because the deleted edge can not be used by any extension.

Suppose now that $G$ does not contain $P_k$ as a subgraph. 
Then every copy of $(H',s,t)$ is simplicial in $G$ since an extension of such a copy gives a subgraph isomorphic to $P_k$.
\end{proof}

We conjecture the following generalization of \cref{prop:isolated-root}.

\begin{conjecture}
    Let $(H,s,t)$ be a two-rooted inherent graph such that $s\neq t$, $d_H(s) = 1$, and let $H'$ be obtained from $H$ by deleting the unique edge incident with $s$. 
    Then $(H',s,t)$ is inherent. 
\end{conjecture}

Note that if $H$ is connected, the fact that $s$ is a pendant vertex is a consequence of the NCI property (recall \cref{obs:NCI-NNI}).

\section{Confining results}\label{sec:confining}

In this section we present two main techniques for constructing confining graphs, which are used to prove that a two-rooted graph is not inherent. 

\subsection{Confining by bipartite graphs}\label{subsec:conf-bipartite}

We first demonstrate how to use complete bipartite graphs to confine trees with imbalanced color classes.
Recall that for any inherent connected two-rooted graph $(H,s,t)$, the graph $H$ is a tree and any tree can be properly colored with two colors.

\begin{thm}\label{thm:Bi-confining}
Let $(H,s,t)$ be a connected inherent two-rooted graph such that 
$H$ is a tree with color classes $C$, $\bar C$, and assume that $s,t\in C$.
Then $|\bar C| \le |C|\le |\bar C|+1$.
\end{thm}

\begin{proof}
Suppose that the chain of inequalities $|\bar C| \le |C|\le |\bar C|+1$ does not hold.
In this case, we show that $(H,s,t)$ is confined by the complete bipartite graph $K_{n,m}$, where $n= |C|$, $m=|\bar C|+2$.
We prove this statement by considering two cases.

\paragraph{Case 1: \boldmath $|C| \geq |\bar C| +2$.}
   Note that for every copy $\cp Hst$ of $(H,s,t)$ in $K_{n,m}$, vertices of the same color class are in the same part of $K_{n,m}$. 
   Hence, $\cpi s$ and $\cpi t$ are in the same part and all vertices of this part are vertices of the copy $\cp Hst$. 
   Therefore we have two free vertices in the other part and the copy can be extended. 
   But this extension can not be closed for obvious reasons. 

\paragraph{Case 2: \boldmath $|C|<|\bar C|$.} 
   Again, for every copy of $H$ in $K_{n,m}$, vertices of the same color class are in the same part of $K_{n,m}$. 
   In particular, $\cpi s$ and $\cpi t$ are in the same part.
   This part should have $n$ vertices, bacause otherwise it should contain all vertices of the color class $\bar C$, but $|\bar C|>n$.
   Therefore, we have two free vertices in the other part and the copy can be extended. 
   But this extension can not be closed for obvious reasons. 
\end{proof}

Let us note that the assumption that $s,t\in C$ in \cref{thm:Bi-confining} is not necessary. 
For example, in an inherent path $T_1(l,0,l')$ where $l\ge 2$ is even and $l'\ge 0$, the distance between the roots is odd, hence, they cannot both belong to $C$.

\subsection{Confining by path augmentation}\label{sec:direct-confinement}

An alternative  way to confine a two-rooted graph is to add paths to it. 
This technique will be used to construct confining graphs for several families, including leaf-extended trees, combs of type II, and rakes.
Before formulating results of this type we need the following additional concept.  

\begin{figure}[!h]
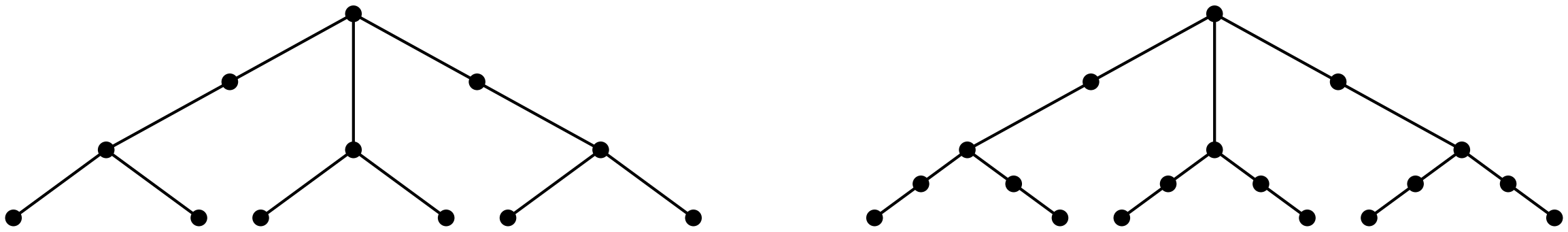

\centerline{\mpfile{pi2rg}{822}\hskip1.6cm \mpfile{pi2rg}{823}}
\caption{A leaf-extended tree (to the right) obtained from a tree (to the left) by subdividing pendant edges, $k=1$.}\label{fig:example}
\end{figure}

\begin{defn}\label{def:leaf-extended}
  A tree is called \emph{leaf-extended} if it can be obtained from a subcubic tree in which every neighbor of every leaf is of degree~$3$, by subdividing each pendant edge $k$ times, for some fixed $k\ge 1$ (see, e.g.,~\cref{fig:example}, where $k=1$). 
\end{defn}

\begin{prop}\label{leaf-extended-trees:confined}
    Let $H$ be a leaf-extended tree, $s$ be a pendant vertex of $H$, and let $t$ be adjacent to a pendant vertex distinct from $s$. 
    Then $(H,s,t)$ can be confined.
\end{prop}

\begin{figure}[!h]
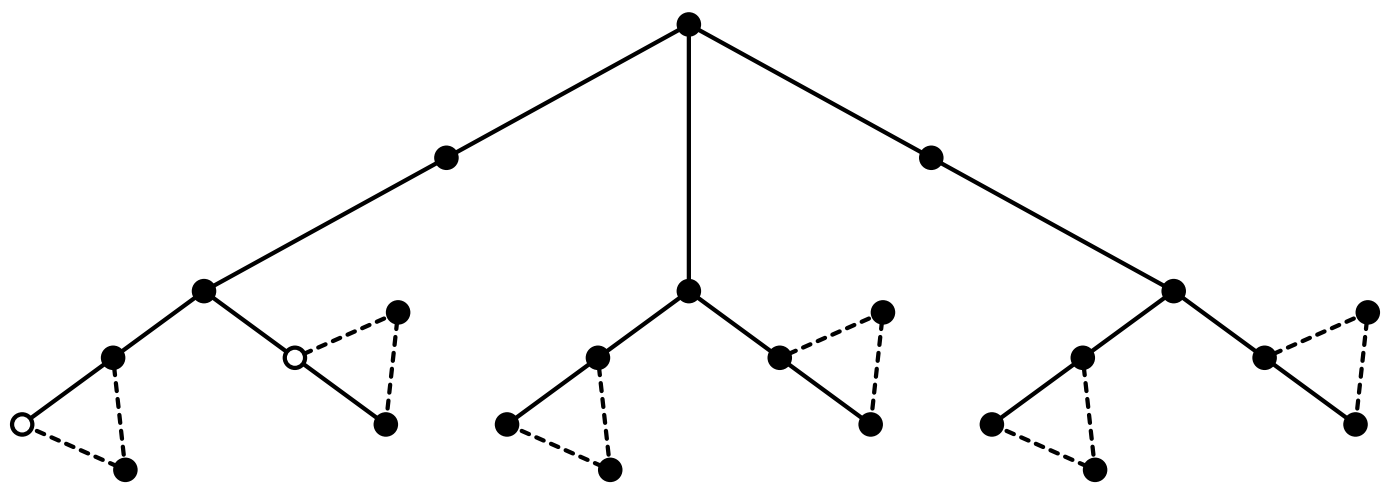

\centerline{\mpfile{pi2rg}{824}}
\caption{A confining graph of the leaf-extended tree from~\cref{fig:example}.
Up to automorphism, the graph admits only one copy of the  leaf-extended  tree, which is not avoidable, since the two roots $s$ and~$t$ (marked by white circles) are not adjacent, thus showing confinement in this case.}\label{fig:trees-confined}
\end{figure}

\begin{proof}
    Let $H$ be a leaf-extended tree obtained from a subcubic tree $T$ in which every neighbor of every leaf is of degree~$3$ by subdividing each pendant edge $k$ times. 
    We construct the confining graph $G$ of $(H,s,t)$ by adding paths of length $2$ through new vertices between each pendant vertex and its neighbor (see \cref{fig:trees-confined}).
    Let $\cp Hst$ be a copy of $(H,s,t)$ in~$G$.
    We claim that from each pair of vertices of degree~$2$ in each triangle of $G$, exactly one of these two vertices belongs to $\cpi H$.
    Indeed, there cannot be more than one such vertex, since each pendant vertex of a leaf-extended tree is at distance $k+1$ from a vertex of degree~$3$.
    Moreover, by a counting argument applied to all such pairs, we infer that there must be exactly one.
    This means that each copy of $(H,s,t)$ has a (unique) extension in $G$, but these extensions are not closable since the roots are in different triangles. 
\end{proof}

Note that for a leaf-extended tree shown in \cref{fig:trees-confined} the distance between the roots is $3$. Thus \cref{thm:Bi-confining} cannot be applied to this two-rooted graph.

\begin{rem}
The proof of \cref{leaf-extended-trees:confined} can be adjusted to the induced case as well. 
However, the confining graphs differ. 
They are obtained from a leaf-extended tree by adding paths of length~$3$ between pendant vertices and their neighbors.  
\end{rem}

In \cite{gurvich2023avoidability} it was shown that, for the induced setting, two-rooted combs of type II, rakes, and split rakes can be confined by adding paths to the graph $H$\footnote{More accurately, in case $p=1$, combs of type II do not satisfy NCI property.}.
A similar construction works also for the NNI setting.

\begin{figure}[!h]
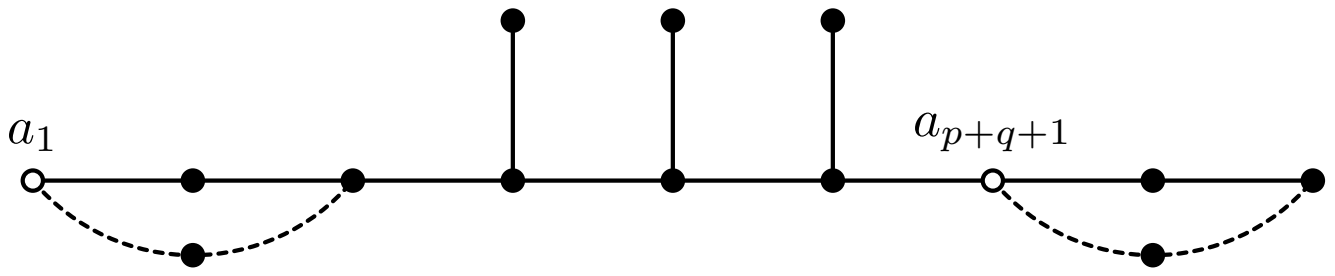

\centerline{\mpfile{pi2rg}{22}}
\caption{Up to automorphism, the graph admits only one copy of $T_2(3,3)$, which is clearly not avoidable, thus showing that $T_2(3,3)$ is confined.}\label{fig:combII-confined}
\end{figure}

\begin{prop}\label{prop:typeII}
    All two-rooted combs of type II can be confined.
\end{prop}
\begin{proof}
Consider $T_2(p,q)$, a two-rooted comb of type II, with $p,q\ge 1$.
If $p=1$, both roots have degree~$1$, so $T_2(p,q)$ does not satisfy the NCI property and is, hence, not inherent, by \Cref{obs:NCI-NNI}.
If $p\ge 2$, add two additional paths of length~$2$ and identify their endpoints with pairs $(a_1,a_{p})$ and $(a_{p+q+1},a_{2p+q})$, respectively (see \Cref{fig:combII-confined}).
\end{proof} 

Note that in some cases \cref{prop:typeII} follows from \cref{thm:Bi-confining}.  
However, for some  combs of type II the assumptions of 
\cref{thm:Bi-confining} fail; in particular, whenever $p+q$ is odd, the distance between roots in $T_2(p,q)$ is odd.

\begin{figure}[!h]
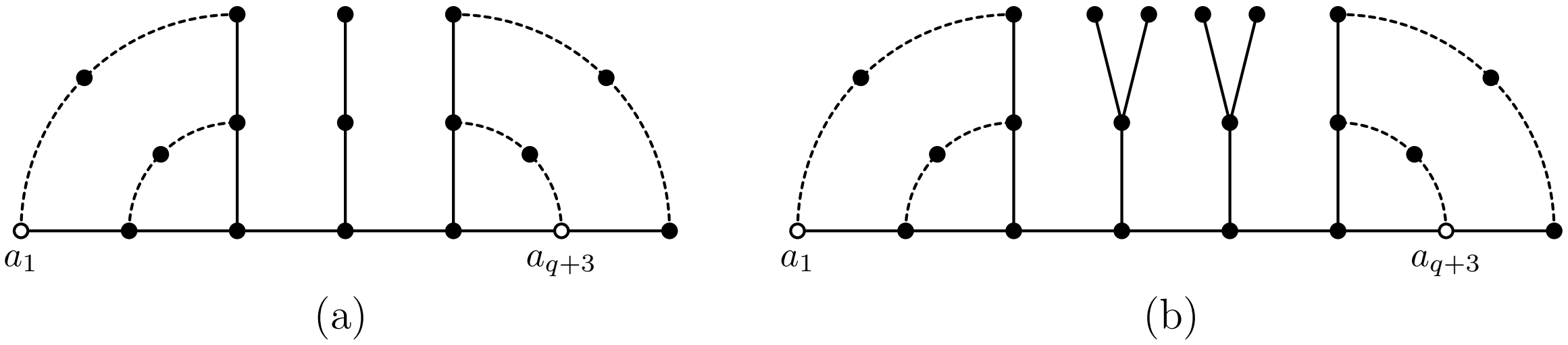

  \centerline{\begin{tabular}{c@{\hskip1cm}c}
    \mpfile{pi2rg}{311} &\mpfile{pi2rg}{722}\\
 (a) &  (b)
    \end{tabular}}
  \caption{Up to automorphism, the graph on the left admits only one copy of $T_4(3)$, while the graph on the right admits only one copy of $T_5(4)$.
  Clearly, neither of these copies is avoidable, hence, the two-rooted graphs $T_4(3)$ and $T_5(3)$ are confined.}\label{fig:rakes-confined}
\end{figure}

Two-rooted rakes and two-rooted split rakes are confined by similar constructions, shown on \cref{fig:rakes-confined}.

Again, in some cases the confining of rakes follows from \cref{thm:Bi-confining}, but in some cases the assumptions of 
\cref{thm:Bi-confining} fail; in particular, whenever $q$ is odd, the distance between roots in $T_4(q)$ and $T_5(q)$ is odd.

\section{Conclusion and Open Questions}\label{sec:conclusion}

In this paper, we extended the study of inherent two-rooted graphs from the induced to the more general not necessarily induced (NNI) setting. 
We developed a new proof technique focused around the CCS (Complementary Connected or Simplicial) property, which is a powerful sufficient condition for inherence. 
Using this framework, we showed that the Necessary Condition for Inherence (NCI) is also sufficient for all two-rooted paths and forks. 
To complement these positive results, we introduced two methods for proving non-inherence: confinement by bipartite graphs with imbalanced partitions and confinement by path augmentation, which we used to rule out several classes of trees.
    
This work opens several avenues for future research. 
In \cite{GKMV20+}, the concept of \textsl{shifting} was used to prove inherence of paths in the induced setting. 
Shifting guarantees that any copy can be transformed into an avoidable one through a series of local modifications.

\begin{question}
    Can our inherence results for paths and forks (\cref{thm:ll',thm:forks}) be strengthened to shifting-type results in the NNI setting?
\end{question}

A central question is the computational complexity of recognizing inherence.
Note that PE-inherence is decidable due to Rabin's theorem about decidability of monadic second-order theory for trees~\cite{rabin1969}.

\begin{question}
    Is the problem of deciding whether a given two-rooted graph is inherent decidable? 
    If so, what is its complexity?
\end{question}

The CCS property is a sufficient condition for inherence. 
It would be interesting to know if it is also necessary.

\begin{question}
    Does there exist an inherent two-rooted tree that is not a CCS graph?
\end{question}

Finally, this work lays the groundwork for a broader classification of inherent graphs.

\begin{question}
    Beyond paths, can we find other families of two-rooted trees whose inherence is characterized by the NCI property? 
\end{question}

\begin{question}
    Is there a more general theory of confinement that applies to a wider range of graph structures?     For example, one can consider more than two roots; for the case of PE-inherence, the first step in this direction was made in \cite{amoeba2024}.
\end{question}

Another open question concerns a relation between the induced and NNI settings. 
It is clear that an induced copy of a two-rooted graph is an NNI copy also, but not necessarily the other way around.
The same holds for extensions. 
Hence, the corresponding two variants of avoidability are incomparable. 
\cref{rem:terminology} shows that some two-rooted graphs are inherent in the NNI setting but not in the induced one.

\begin{question}
    Is there a two-rooted graph that is inherent in the induced setting but is not inherent in the NNI setting?
\end{question}

{\huge

\bigskip

\newpage
\appendix
\noindent Appendix

}

\section{Dictionary of abbreviations}\label{sec:dictionary}
  
For convenience of the reader, we collect here a list of abbreviations of various properties used in this paper.
All these properties apply to a two-rooted graph $\cH =(H,s,t)$, except for the $H$-CC property, which depends on a graph $H$ and applies to a graph $G$.

\begin{description}
\item [AS:] Arc Separating.
In any $\cH$-smooth graph $G=(V, E)$ and any arc $e=(x,y)\in E$ such that $e$ is contained in a copy of $H$ in $G$ there exists a copy of $H$ that contains $x$ and does not contain $y$ (see \cref{def:H-universality}).
\item [AU:] Arc Universal. 
In any $\cH$-smooth graph $G$, any arc of $G$ plays a role of each arc of $H$ (see \cref{def:H-universality}). 
\item [CCS:] every $\cH$-smooth  graph has the $H$-CC property (see \cref{def:CCS}).
\item [EM:] Extension Movable. 
The extension graph of $\cH$ contains more than one copy of $H$ (see \cref{def:EM}). 
\item [{\boldmath{$H$-CC property}} {\mdseries(of a graph $G$)}:] there exists a copy $H'$ of $H$ in $G$ such that $G-V(H')$ is connected (see  \cref{def:HCC}).
\item [HCS:] H-Component is Small. See \cref{def:HCS}.
\item [NCI:] Necessary Condition for Inherence. 
The graph $H$ is either a path or a tree with maximum degree~$3$ such that $s$ is a pendant vertex and $t$ either coincides with $s$ or has degree $2$ (see \cref{def:NCI}). 
\item [VU:] Vertex Universal. 
In any $\cH$-smooth graph $G$, any vertex of $G$ plays a role of each vertex of $H$ (see \cref{def:H-universality}). 
\item [WAU:] Weakly Arc Universal. 
In any $\cH$-smooth graph $G$, any arc of $G$ contained in a copy of $H$ in $G$ plays a role of each arc of $H$ (see \cref{def:H-universality}). 
\item [WVU:] Weakly Vertex Universal. 
In any $\cH$-smooth graph $G$, any vertex of $G$ contained in a copy of $H$ in $G$ plays a role of each vertex of $H$ (see \cref{def:H-universality}). 
\end{description}

\newpage
\section{Two-rooted graphs with \boldmath $n=4,5,6$ vertices}\label{sec:small-graphs}

As shown in \cref{thm:ll'}, all NCI two-rooted paths are inherent.
Here we list small trees that are not covered by these general results. 
In particular, $H$ should have vertices of degree~$3$.

\subsection{\boldmath $n=4$}

The only $4$-vertex tree that is not a path is the star $K_{1,3}$. 
The only possible placement of root for an inherent $K_{1,3}$ satisfying the NCI property is a double root in a pendant vertex. 
But this two-rooted graph is not inherent. 
It is confined by a bipartite graph. 
In terms of \cref{thm:Bi-confining} $|C|=3$, $|\bar C|=1< |C|-1$. 
Hence, the graph is not inherent due to \cref{thm:Bi-confining}.

\subsection{\boldmath $n=5$}\label{subsec:forks}

Each subcubic tree on $5$~vertices that is not a path can be obtained from a subcubic tree on $4$~vertices by attaching a pendant edge to a vertex of degree at most $2$.
There exists a unique such tree on $5$~vertices, namely, the fork graph.

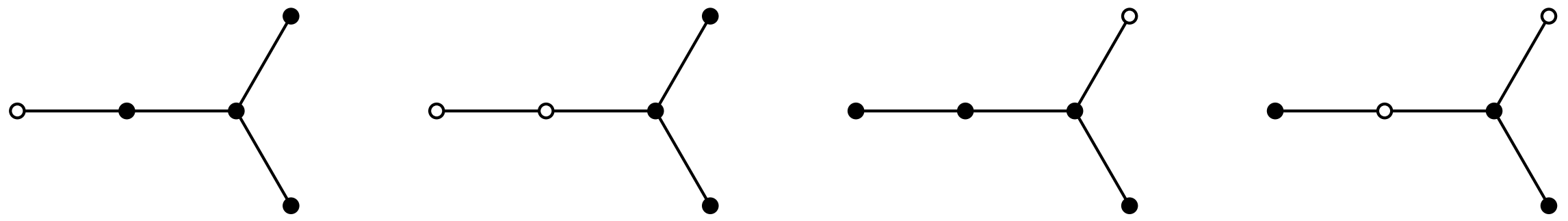
\begin{figure}[!h]

  \centering

  \begin{tabular}{p{0.22\textwidth}p{0.22\textwidth}p{0.22\textwidth}p{0.22\textwidth}}
\mpfile{trees2rooted}{50}&\mpfile{trees2rooted}{51}&\mpfile{trees2rooted}{52}&\mpfile{trees2rooted}{53}\\[6pt]
not inherent &\multicolumn{3}{c}{all these are PE-forks and they are inherent by \cref{thm:forks}}
  \end{tabular}
  \caption{All the $5$-vertex two-rooted graphs satisfying the NCI property. 
  Roots are shown as white circles. 
If only one root is shown, it means that the root is double, that is, $s=t$.}\label{pic:2rooted-forks}
\end{figure}

To determine the number of NCI two-rooted forks up to isomorphism, note that there are two~orbits of fork automorphisms on pendant vertices and only one vertex of degree $2$. 
Thus, the number of two-rooted graphs in the list is $2+2=4$.  
The corresponding two-rooted forks are shown in \cref{pic:2rooted-forks}. 

\begin{figure}[!h]
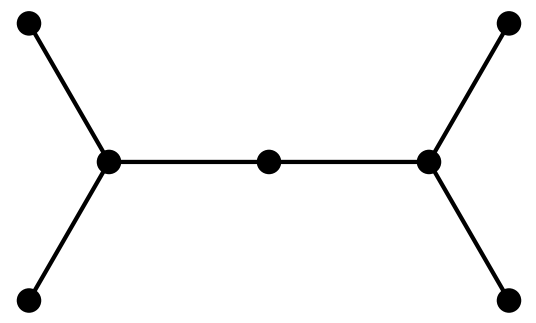

  \centering
  \mpfile{trees2rooted}{1}
  \caption{Confining graph for the non-inherent two-rooted fork in \cref{pic:2rooted-forks}.}\label{pic:cyrzh}
\end{figure}

The first two-rooted fork in \cref{pic:2rooted-forks} is confined by the graph shown in \cref{pic:cyrzh}. 
Hence, it is not PE-inherent. 
Up to automorphism, there is only one copy of the two-rooted fork in the graph, and this copy is not avoidable.  

\subsection{\boldmath $n=6$}

\begin{figure}[!h]
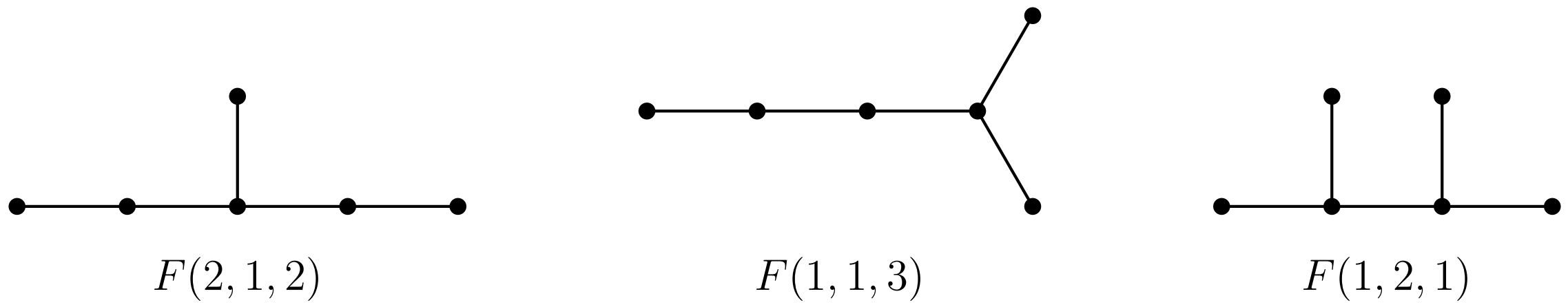

  \centering
  \begin{tabular}{c@{\qquad\qquad}c@{\qquad\qquad}c}
\mpfile{trees2rooted}{61}&  \mpfile{trees2rooted}{62}&  \mpfile{trees2rooted}{63}\\[6pt]
$F(2,1,2)$& $F(1,1,3)$& $F(1,2,1)$
  \end{tabular}
  \caption{$6$-vertex non-path trees that up to isomorphism.}\label{pic:6trees}
\end{figure}

Each subcubic tree on $6$~vertices that is not a path can be obtained from a subcubic tree on $5$~vertices  by attaching a pendant edge to a vertex of degree at most $2$.
Up to isomorphism, this results in three trees, see \cref{pic:6trees}. 
We provide their names in our tree naming system. 

\begin{figure}[!h]
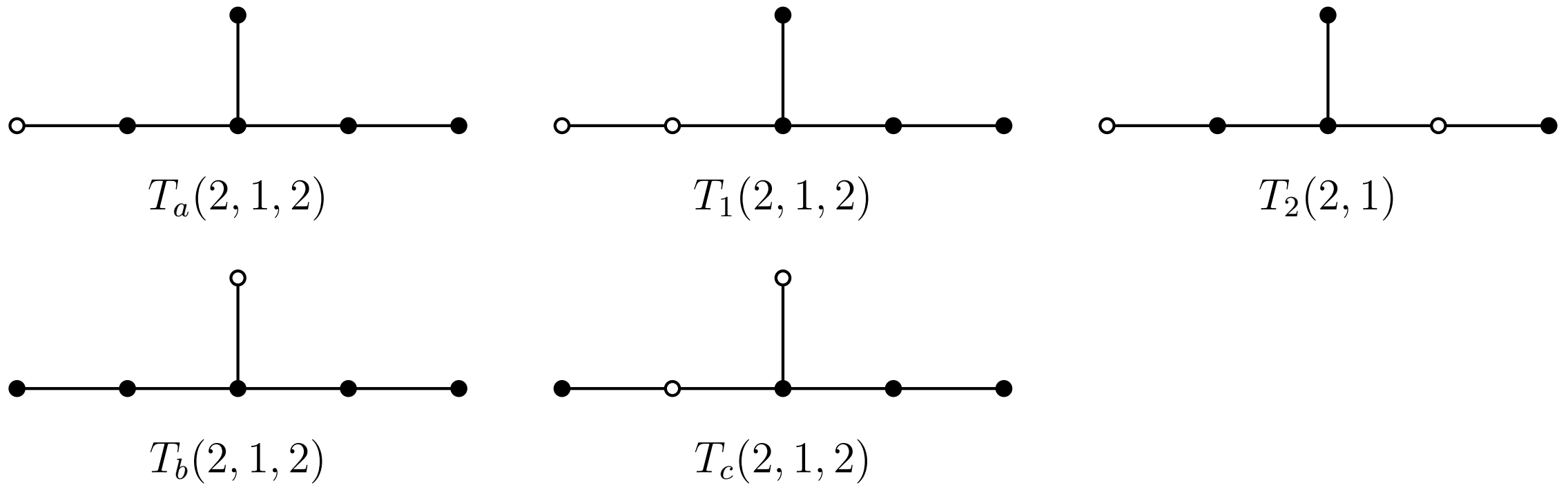


  \centering

  \begin{tabular}{c@{\qquad}c@{\qquad}c}
\mpfile{trees2rooted}{610}&\mpfile{trees2rooted}{611}&\mpfile{trees2rooted}{612}\\[6pt]
$T_a(2,1,2)$ & $T_1(2,1,2)$& $T_2(2,1)$\\[12pt]
\mpfile{trees2rooted}{613}&\mpfile{trees2rooted}{614}&\\[6pt]
$T_b(2,1,2)$ & $T_c(2,1,2)$&\\[12pt]
  \end{tabular}
  \caption{Two-rooted graphs $(H,s,t)$ with $H=F(2,1,2)$ satisfying the NCI property.
  Among them, only $T_1(2,1,2)$ is inherent.}\label{pic:F212}
\end{figure}

 There are five two-rooted graphs on the graph $F(2,1,2)$ satisfying the NCI property. 
 Indeed, there are two different pendant vertices up to automorphism and two vertices of degree~$2$. 
 For an endpoint of path $P_5$ the two vertices of degree $2$ differ, while for a third pendant vertex they are the same up to automorphism of one-rooted graph.   
These two-rooted graphs, along with their names, are shown in \cref{pic:F212}.

For a proof that $T_1(2,1,2)$ is inherent, see \cref{sec:T1212}. 
From \cref{prop:typeII} we conclude that $T_2(2,1)$ is not inherent. 
The leftmost two-rooted graph in the first row of \cref{pic:F212} is non-standard for our naming conventions; we denote it as $T_a(2,1,2)$. 
It is confined by the graph shown in \cref{pic:6confinement}(a). 
In this graph, the only vertex of degree $3$ that is an endpoint of two edge-disjoint paths of length $2$ is the central vertex. 
Hence, there are exactly two copies of  $T_a(2,1,2)$ in the graph and neither of them is avoidable.

The two graphs in the second row of \cref{pic:F212} are not inherent.
They are both confined by the graph shown in \cref{pic:6confinement}(b). 
In this graph, there is a unique path of length $4$ up to automorphism. 
Consequently, up to automorphism, there is only one copy of $T_b(2,1,2)$ in the graph, and this copy is not avoidable. 
Similar observations apply to $T_c(2,1,2)$.

\begin{figure}[!h]
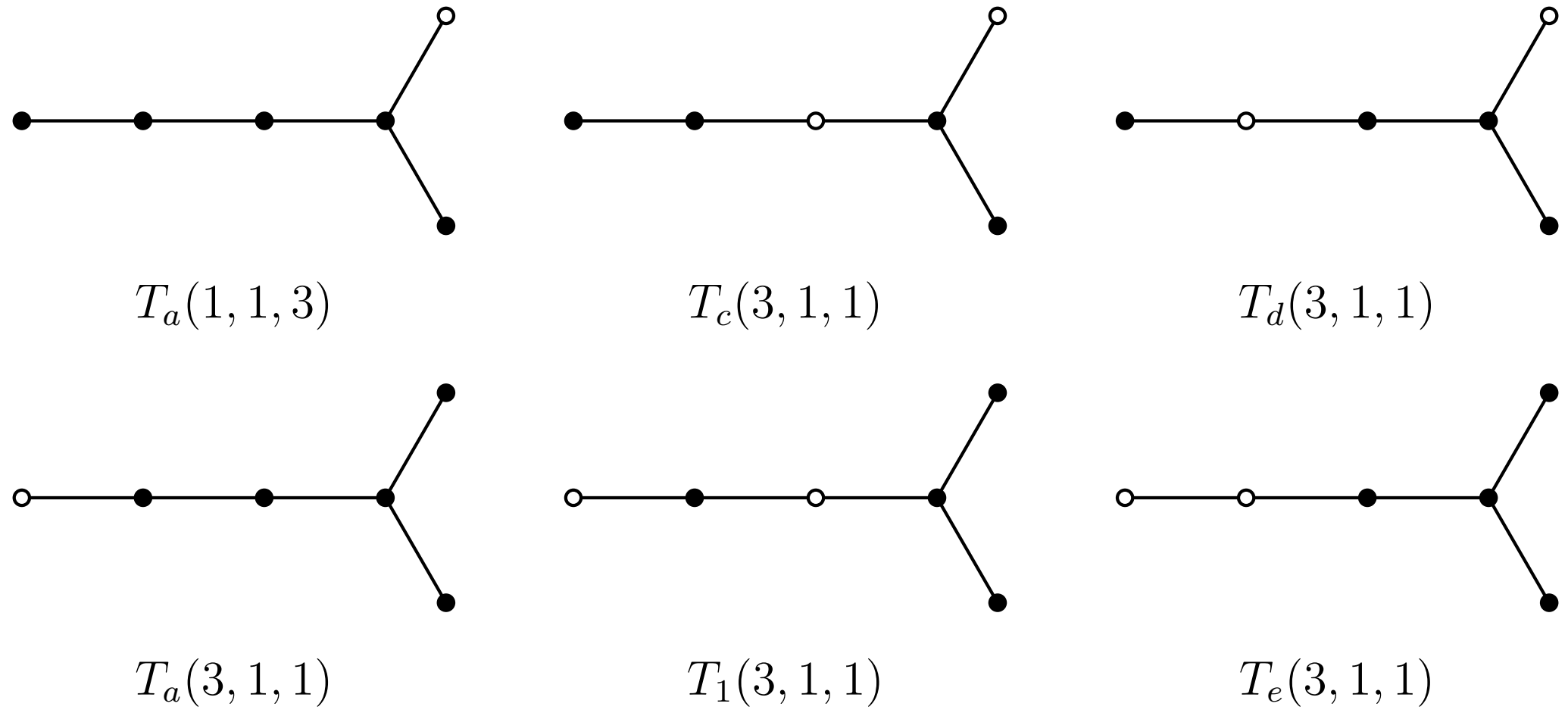

  \centering
  \begin{tabular}{c@{\hspace*{0.95cm}}c@{\hspace*{0.95cm}}c}
\mpfile{trees2rooted}{620}&\mpfile{trees2rooted}{622}&\mpfile{trees2rooted}{621}\\[6pt]
$T_a(1,1,3)$ & $T_c(3,1,1)$ & $T_d(3,1,1)$ \\[12pt]
\mpfile{trees2rooted}{623}&\mpfile{trees2rooted}{625}&\mpfile{trees2rooted}{624}\\[6pt]
$T_a(3,1,1)$  & $T_1(3,1,1)$ &$T_e(3,1,1)$ \\[12pt] 
  \end{tabular}
  \caption{Two-rooted graphs $(H,s,t)$ with $H=F(1,1,3)$ satisfying the necessary condition. 
  None of these two-rooted graphs is inherent.}\label{pic:F113}
\end{figure}

There are six two-rooted graphs on the graph $F(1,1,3)$ satisfying the NCI property. 
Indeed, up to automorphism, there are two different pendant vertices and two vertices of degree~$2$. 
Each pendant vertex gives three variants. 
All variants are shown in \cref{pic:F113}. 
It turns out that none of these two-rooted graphs is inherent.

Graphs $T_a(1,1,3)$, $T_a(3,1,1)$, $T_c(3,1,1)$, and $T_1(3,1,1)$ are not inherent, since \cref{thm:Bi-confining} is applied for them: distance between roots is even. 
For these graphs we have $|C|=4$, $|\bar C|=2$, and inequalities in  \cref{thm:Bi-confining} are violated. 
A graph confining $T_d(3,1,1)$ and $T_e(3,1,1)$ is shown in \cref{pic:6confinement}(c).
Up to automorphism, there is a unique copy of $T_d(3,1,1)$ in this graph and that copy is not avoidable. 
The same holds for $T_e(3,1,1)$.

\begin{figure}[!h]
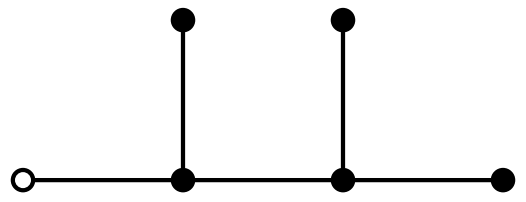

  \centering
\mpfile{trees2rooted}{631}
  \caption{The only two-rooted graph $(H,s,t)$ with $H=F(1,2,1)$  satisfying the necessary condition is $T_1(1,2,1)$. 
  It is not inherent.}\label{pic:F121}
\end{figure}

Since $F(1,2,1)$ has no vertices of degree $2$ and its pendant vertices form an orbit of the automorphism group, there is exactly one two-rooted graph on the graph $F(1,2,1)$ satisfying the necessary condition. 
It is $T_1(1,2,1)$ (see \Cref{pic:F121}). 
This two-rooted graph is not inherent, since it is confined by Petersen graph, see \cref{pic:6confinement}(d). Similarly to most cases considered before, up to automorphism, there is only one copy of $T_1(1,2,1)$.
It is easy to verify using   \cref{pic:6confinement}(d) that this copy is not avoidable.

\begin{figure}
  \centering

  \begin{tabular}{c@{\qquad}c@{\qquad}c}
   \raisebox{25pt}{\mpfile{trees2rooted}{11}}&\mpfile{trees2rooted}{13}\\[6pt]
(a) Confining $T_a(2,1,2)$  & (b) Confining $T_b(2,1,2)$ and $T_c(2,1,2)$\\[12pt]
\raisebox{14pt}{\mpfile{trees2rooted}{22}}&
\mpfile{trees2rooted}{2}   \\[6pt]
 (c) Confining $T_d(3,1,1)$ and $T_e(3,1,1)$ & (d) Petersen confines $T_1(1,2,1)$ 
  \end{tabular}
  \caption{Confining graphs for the two-rooted graphs mentioned in the text.}\label{pic:6confinement}
\end{figure}

\bibliographystyle{plain}
\bibliography{shifting}

@article {GKMV20+,
    AUTHOR = {Gurvich, Vladimir and Krnc, Matja\v{z} and Milani\v{c}, Martin and
              Vyalyi, Mikhail},
     TITLE = {Shifting paths to avoidable ones},
   JOURNAL = {J. Graph Theory},
  FJOURNAL = {Journal of Graph Theory},
    VOLUME = {100},
      YEAR = {2022},
    NUMBER = {1},
     PAGES = {69--83},
     DOI = {10.1002/jgt.22766}, 
}

@article{MR2204109,
	Author = {Heggernes, Pinar},
	Doi = {10.1016/j.disc.2005.12.003},
	Fjournal = {Discrete Mathematics},
	Journal = {Discrete Math.},
	Mrclass = {05C85 (05C10 05C75)},
	Mrnumber = {2204109},
	Mrreviewer = {Dieter Kratsch},
	Number = {3},
	Pages = {297--317},
	Title = {Minimal triangulations of graphs: a survey},
	Url = {https://doi.org/10.1016/j.disc.2005.12.003},
	Volume = {306},
	Year = {2006}}

@article{MR2057267,
	Author = {Berry, Anne and Blair, Jean R. S. and Heggernes, Pinar and Peyton, Barry W.},
	Doi = {10.1007/s00453-004-1084-3},
	Fjournal = {Algorithmica. An International Journal in Computer Science},
	Journal = {Algorithmica},
	Mrclass = {68R10 (05C85 68W05)},
	Mrnumber = {2057267},
	Number = {4},
	Pages = {287--298},
	Title = {Maximum cardinality search for computing minimal triangulations of graphs},
	Url = {https://doi.org/10.1007/s00453-004-1084-3},
	Volume = {39},
	Year = {2004}}

@article{MR1927566,
	Author = {Chv\'{a}tal, Va\v{s}ek and Rusu, Irena and Sritharan, R.},
	Doi = {10.1016/S0012-365X(01)00166-2},
	Fjournal = {Discrete Mathematics},
	Journal = {Discrete Math.},
	Mrclass = {05C75},
	Mrnumber = {1927566},
	Number = {1-2},
	Pages = {445--448},
	Title = {Dirac-type characterizations of graphs without long chordless cycles},
	Url = {https://doi.org/10.1016/S0012-365X(01)00166-2},
	Volume = {256},
	Year = {2002}}

@article{MR485552,
	Author = {Ohtsuki, Tatsuo and Cheung, Lap Kit and Fujisawa, Toshio},
	Doi = {10.1016/0022-247X(76)90182-7},
	Fjournal = {Journal of Mathematical Analysis and Applications},
	Journal = {J. Math. Anal. Appl.},
	Mrclass = {05C35},
	Mrnumber = {485552},
	Number = {3},
	Pages = {622--633},
	Title = {Minimal triangulation of a graph and optimal pivoting order in a sparse matrix},
	Url = {https://doi.org/10.1016/0022-247X(76)90182-7},
	Volume = {54},
	Year = {1976}}

@article{MR130190,
	Author = {Dirac, G. A.},
	Doi = {10.1007/BF02992776},
	Fjournal = {Abhandlungen aus dem Mathematischen Seminar der Universit\"{a}t Hamburg},
	Journal = {Abh. Math. Sem. Univ. Hamburg},
	Mrclass = {05.40},
	Mrnumber = {130190},
	Mrreviewer = {F. Harary},
	Pages = {71--76},
	Title = {On rigid circuit graphs},
	DOI = {10.1007/BF02992776},
	Volume = {25},
	Year = {1961},
	}

@article{MR4245221,
	Author = {Bonamy, Marthe and Defrain, Oscar and Hatzel, Meike and Thiebaut, Jocelyn},
	Doi = {10.37236/9030},
	Fjournal = {Electronic Journal of Combinatorics},
	Journal = {Electron. J. Combin.},
	Mrclass = {05C38 (05C60)},
	Mrnumber = {4245221},
	Number = {4},
	Pages = {Paper No. 4.46, 9},
	Title = {Avoidable paths in graphs},
	Url = {https://doi.org/10.37236/9030},
	Volume = {27},
	Year = {2020}}

@inproceedings{BCGMS19,
  author       = {Jesse Beisegel and
                  Maria Chudnovsky and
                  Vladimir Gurvich and
                  Martin Milani{\v c} and
                  Mary Servatius},
  editor       = {Zachary Friggstad and
                  J{\"{o}}rg{-}R{\"{u}}diger Sack and
                  Mohammad R. Salavatipour},
  title        = {Avoidable Vertices and Edges in Graphs},
  booktitle    = {Algorithms and Data Structures - 16th International Symposium, {WADS}
                  2019, Edmonton, AB, Canada, August 5-7, 2019, Proceedings},
  series       = {Lecture Notes in Computer Science},
  volume       = {11646},
  pages        = {126--139},
  publisher    = {Springer},
  year         = {2019},
  doi          = {10.1007/978-3-030-24766-9\_10},
}

@article{BCGMS19DAM,
	Author = {Beisegel, Jesse and Chudnovsky, Maria and Gurvich, Vladimir and Milani\v{c}, Martin and Servatius, Mary},
	Fjournal = {Discrete Applied Mathematics. The Journal of Combinatorial Algorithms, Informatics and Computational Sciences},
	Journal = {Discrete Appl. Math.},
	Pages = {285--300},
	Title = {Avoidable vertices and edges in graphs: {E}xistence, characterization, and applications},
	Volume = {309},
	Year = {2022},
     DOI = {10.1016/j.dam.2021.12.006},
    }

@article {MR2191642,
    AUTHOR = {Berry, Anne and Bordat, Jean-Paul and Heggernes, Pinar and
              Simonet, Genevi\`eve and Villanger, Yngve},
     TITLE = {A wide-range algorithm for minimal triangulation from an
              arbitrary ordering},
   JOURNAL = {J. Algorithms},
  FJOURNAL = {Journal of Algorithms. Cognition, Informatics and Logic},
    VOLUME = {58},
      YEAR = {2006},
    NUMBER = {1},
     PAGES = {33--66},
       DOI = {10.1016/j.jalgor.2004.07.001},
       URL = {https://doi.org/10.1016/j.jalgor.2004.07.001},
}

@article{amoeba2024,
    AUTHOR = {Gurvich, Vladimir and Krnc, Matja\v{z} and Vyalyi, Mikhail},
     TITLE = {Growing trees and amoebas' replications},
   JOURNAL = {Results Math.},
  FJOURNAL = {Results in Mathematics},
    VOLUME = {80},
      YEAR = {2025},
    NUMBER = {4},
     PAGES = {Paper No. 117, 17},
       DOI = {10.1007/s00025-025-02421-6},
}

@article {gurvich2023avoidability,
    AUTHOR = {Gurvich, Vladimir and Krnc, Matja{\v z} and Milani{\v c}, Martin
              and Vyalyi, Mikhail},
     TITLE = {Avoidability beyond paths},
   JOURNAL = {Electron. J. Combin.},
  FJOURNAL = {Electronic Journal of Combinatorics},
    VOLUME = {32},
      YEAR = {2025},
    NUMBER = {2},
     PAGES = {Paper No. 2.27, 39},
      MRCLASS = {05C75 (05C38 05C60 05C63)},
  MRNUMBER = {4904889},
       DOI = {10.37236/13023},
}

@book {diestel-graph-theory,
    AUTHOR = {Diestel, Reinhard},
     TITLE = {Graph theory},
    SERIES = {Graduate Texts in Mathematics},
    VOLUME = {173},
   EDITION = {Sixth},
 PUBLISHER = {Springer, Berlin},
      YEAR = {2025},
   MRCLASS = {05-01 (05Cxx)},
  MRNUMBER = {4874150},
}

@article{rabin1969,
    AUTHOR = {Rabin, Michael O.},
     TITLE = {Decidability of second-order theories and automata on infinite trees},
   JOURNAL = {Trans. Amer. Math. Soc.},
  FJOURNAL = {Transactions of the American Mathematical Society},
    VOLUME = {141},
      YEAR = {1969},
     PAGES = {1--35},
       DOI = {10.1090/S0002-9947-1969-0246760-1},
}
\end{document}